# Implicit Functions from Topological Vector Spaces to Banach Spaces

Helge Glöckner

**Abstract.** We prove implicit function theorems for mappings on topological vector spaces over valued fields. In the real and complex cases, we obtain implicit function theorems for mappings from arbitrary (not necessarily locally convex) topological vector spaces to Banach spaces.

## Introduction

In this article, we prove implicit function theorems (and generalizations) for mappings from topological vector spaces over valued fields to Banach spaces. Our main results can be summarized as follows. Let $(\mathbb{K}, |.|)$ be a (non-discrete, not necessarily complete) valued field, $E$ and $F$ be topological $\mathbb{K}$-vector spaces, $U \subseteq E$ and $V \subseteq F$ be open subsets, and $f : U \times V \to F$ be a mapping which is of class $C^1$. Let $(x_0, y_0) \in U \times V$ such that $d_2 f(x_0, y_0, \bullet) \in \mathrm{GL}(F)$. Then, given the respective hypotheses stated in the first four columns of the following table, there exists an open neighbourhood $Q \subseteq U$ of $x_0$ and an open neighbourhood $B \subseteq V$ of $y_0$ such that, for every $x \in Q$, there is a unique element $\beta(x) \in B$ such that $f(x, \beta(x)) = f(x_0, y_0)$, and the mapping $\beta \colon Q \to B$ so obtained has the property shown in the last column:

| $\mathbb{K}$ | $E$ | $F$ | $f$ | $\beta$ |
|---|---|---|---|---|
| $\mathbb{R}$ | arbitrary | $\dim F < \infty$ | $C^k$ | $C^k$ |
| $\mathbb{R}$ | arb. | Banach | $C^{k+1}$ | $C^k$ |
| $\mathbb{R}$ | arb. | Banach | $\mathbb{R}$-analytic | $\mathbb{R}$-analytic |
| $\mathbb{C}$ | arb. | Banach | $\mathbb{C}$-analytic | $\mathbb{C}$-analytic |
| local | arb. | $\dim F < \infty$ | $C^1$ | $C^0$ |
| local | metrizable | $\dim F < \infty$ | $C^k$ | $C^k$ |
| complete ultram. | metrizable | Banach | $C^{k+1}$ | $C^k$ |
| $\mathbb{K} \subseteq \mathbb{R}$ | metrizable | Banach | $C^2$ | $C^1$ |
| arb. | arb. | Banach | $C^2$ | $C^0$ |
| arb. | Banach | Banach | $SC^1$ | $SC^1$ |
| complete | Banach | Banach | $SC^k$ | $SC^k$ |

Here $k \in \mathbb{N} \cup \{\infty\}$, and $C^k$-maps are understood in the sense of [2], where a differential calculus over arbitrary non-discrete topological fields is developed. A map between open subsets of real locally convex spaces is $C^k$ in the former sense if and only if it is $C^k$ in the sense of Michal-Bastiani (*i.e.*, a Keller $C_c^k$-map [19]). The symbol $SC^k$ refers to $k$ times strictly differentiable mappings, as defined below.

Our results were inspired by Hiltunen's implicit function theorems, which he formulated in the setting of $C_\Pi^k$-maps on locally convex spaces (see [16] for the real case, [17] for the complex case). In contrast to that paper, we are working throughout in the realm of topology, no recourse to convergence structures is necessary. Furthermore, we are able to work over valued fields other than $\mathbb{R}$ or $\mathbb{C}$, and need not assume that the domains be locally convex.

Our approach is based on the classical idea that every implicit function theorem has an underlying "inverse function theorem with parameters" (cf. [18], [20, Thm. 3.2.1], [34]). For example, in the case of a $C^k$-map $f: U \times V \to F$, where $U$ is a subset of a real topological vector space and $F = \mathbb{R}^n$, given $(x_0, y_0) \in U \times V$ with $d_2 f(x_0, y_0, \bullet)$ invertible, after shrinking $U$ we interpret $f$ as a family $(f_x)_{x \in U}$ of mappings $f_x := f(x, \bullet) \colon V \to F$ between open subsets of the finite-dimensional space $F$, to which the classical inverse function theorem applies, and then show that $\psi \colon (x, y) \mapsto f_x^{-1}(y)$ makes sense on some open neighbourhood of $(x_0, f(x_0, y_0))$ and defines a $C^k$-map there. Then $y(x) = \psi(x, f(x_0, y_0))$ gives a $C^k$-solution to the implicit equation $f(x, y) = f(x_0, y_0)$.

Local convexity does not play a role in this line of argument, and we are also able to tackle the ultrametric case. This is much more difficult, since the absence of a fundamental theorem of calculus and mean value theorem in this setting makes it necessary to discuss continuous extensions to iterated difference quotient maps, rather than the mere existence and continuity of higher differentials. To get from the Banach case to implicit functions on open subsets of metrizable topological vector spaces over a complete ultrametric field $\mathbb{K}$, we exploit the fact that a mapping on a metrizable space is of class $C^k$ if and only if all of its compositions with smooth maps from $\mathbb{K}^{k+1}$ to the space are of class $C^k$ [2]. This result can be seen as an adaptation of ideas from the convenient differential calculus of Frölicher, Kriegl and Michor ([9], [22], in particular also [21]) and Souriau's theory of diffeological spaces ([36], [25]) to non-archimedian analysis.

The ultrametric implicit function theorem developed here makes it child's play to turn the diffeomorphism groups of $\sigma$-compact, finite-dimensional smooth manifolds over totally disconnected local fields into Lie groups (see [12]). For related studies, cf. [26], [27].

Generalizations of the classical inverse function theorem (or implicit function theorem) for special classes of mappings taking values in non-Banach spaces can be found in [14] (for tame Fréchet spaces), [24] (for Silva spaces), and [28] (for so-called "coordinate spaces").

The article is structured as follows. Having described the precise setting of differential calculus used in the present article (Section 1), we present our results in the real and complex case. We do not try to be self-contained here, but rather re-use the standard Inverse

Mapping Theorem for real Banach spaces and its corollaries (as proved in [23]), which should be well-known to most readers, to get to the point as quickly as possible.

In Section 3 we recall the notion of a strictly differentiable mapping from an open subset of a normed vector space over a valued field $\mathbb{K}$ to a polynormed $\mathbb{K}$-vector space (cf. [5]). We show that any strictly differentiable map is of class $C^1$, and we show that every $C^2$-map from an open subset of a normed $\mathbb{K}$-vector space to a polynormed $\mathbb{K}$-vector space is strictly differentiable. In Section 4, we specialize to locally compact $\mathbb{K}$. In this case, a mapping from an open subset of a finite-dimensional $\mathbb{K}$-vector space to a polynormed $\mathbb{K}$-vector space is $C^1$ if and only if it is strictly differentiable, if and only if it is "locally uniformly differentiable" (an *a priori* even stronger differentiability property). In Section 5, we introduce the class of $k$ times strictly differentiable mappings ($SC^k$-maps, for short). Any such map is $C^k$, and, conversely, we show that every $C^{k+1}$-map from an open subset of a normed vector space over a valued field $\mathbb{K}$ to a polynormed $\mathbb{K}$-vector space is $SC^k$.

We recall from [32, Example 26.6] that there exists a function $f\colon \mathbb{Z}_p \to \mathbb{Q}_p$ from the $p$-adic integers to the $p$-adic numbers which is differentiable (in the naïve sense) at each $x \in \mathbb{Z}_p$, with $f'(x) = 1$ (whence $f'(x)$ is invertible, and $f'\colon \mathbb{Z}_p \to \mathbb{Q}_p$ a continuous map), but such that $f$ is not injective on any zero-neighbourhood (which refutes earlier claims in [31]). Thus, in the ultrametric case, an inverse function theorem cannot be based on the mere existence and continuity of differentials. In contrast, the $SC^k$-property is well-adapted to inverse and implicit function theorems. An inverse function theorem for once strictly differentiable mappings between open subsets of Banach spaces over complete valued fields is well-known (see [5, 1.5.1], where no proofs are given and where all Banach spaces over ultrametric fields are assumed ultrametric). In the finite-dimensional ultrametric case, higher order differentiability has been discussed in [1] for implicit functions, in [35] (with merely partial proofs) for inverse functions. Going beyond these known facts, using an inductive argument which goes back and forth between inverse functions and implicit functions, we establish the Inverse Function Theorem and the Implicit Function Theorem for $SC^k$-maps between open subsets of Banach spaces over complete valued fields (Section 7). Combining these results with parameter-dependent Newton approximation (from Section 6) and the specific tools of differential calculus on metrizable spaces outlined above, we obtain a general Implicit Function Theorem for mappings from metrizable topological vector spaces to (not necessarily ultrametric) Banach spaces over complete ultrametric fields (Section 8).

In an appendix, we show that, in the real case, every $k$ times continuously Fréchet differentiable mapping is an $SC^k$-map. For $k = 1$, the converse also holds [5, 2.3.3].

All results are formulated in a way which transports as much useful information as possible. For example, instead of formulating mere implicit function theorems, in accordance with our general philosophy we explicitly spell out "inverse function theorems with parameters," and we provide information concerning the size of images of balls. Such refined information is useful for the discussion of diffeomorphism groups, and also in other contexts (*e.g.*, [13]).

# 1 Differential calculus over topological fields

In this article, we are working in the setting of differential calculus over non-discrete topological fields developed in [2]. In this section, we briefly recall basic definitions and facts.

Unless stated otherwise, in this section $\mathbb{K}$ denotes a non-discrete topological field. All topological vector spaces are assumed Hausdorff. Before we define $C^k$-maps, we need an efficient notation for the domains of certain mappings associated with $C^k$-maps.

**Definition 1.1** If $E$ is a topological $\mathbb{K}$-vector space and $U \subseteq E$ an open subset, we define $U^{[0]} := U$ and
$$U^{[1]} := \{(x, y, t) \in U \times E \times \mathbb{K} \colon x + ty \in U\},$$
which is an open subset of the topological $\mathbb{K}$-vector space $E \times E \times \mathbb{K}$. Having defined $U^{[j]}$ inductively for a natural number $j \geq 1$, we set $U^{[j+1]} := (U^{[j]})^{[1]}$.

In particular, $E^{[1]} = E \times E \times \mathbb{K}$, $E^{[2]} = E \times E \times \mathbb{K} \times E \times E \times \mathbb{K} \times \mathbb{K}$, etc.

**Definition 1.2** Let $E$ and $F$ be topological $\mathbb{K}$-vector spaces, and $f \colon U \to F$ be a mapping, defined on an open subset $U \subseteq E$. We say that $f$ is *of class* $C_{\mathbb{K}}^0$ if $f$ is continuous, we set $f^{[0]} := f$ and call $f^{[0]}$ the 0-th extended difference quotient map of $f$. If $f$ is continuous and there exists a continuous mapping $f^{[1]} \colon U^{[1]} \to F$ such that

$$\tfrac{1}{t}(f(x+ty) - f(x)) = f^{[1]}(x, y, t) \quad \text{for all } (x, y, t) \in U^{[1]} \text{ such that } t \neq 0, \tag{1}$$

we say that $f$ is *of class* $C_{\mathbb{K}}^1$, and call $f^{[1]}$ the (first) *extended difference quotient map* of $f$. Here $f^{[1]}$ is uniquely determined, as $\mathbb{K}$ is non-discrete. Recursively, having defined $C_{\mathbb{K}}^j$-maps and $j$-th extended difference quotient maps for $j = 0, \ldots, k-1$ for some natural number $k \geq 2$, we call $f$ a mapping *of class* $C_{\mathbb{K}}^k$ if $f$ is of class $C_{\mathbb{K}}^{k-1}$ and $f^{[k-1]}$ is of class $C_{\mathbb{K}}^1$. In this case, we define the $k$-th extended difference quotient map of $f$ via
$$f^{[k]} := (f^{[k-1]})^{[1]} \colon U^{[k]} \to F.$$

The mapping $f$ is *of class* $C_{\mathbb{K}}^\infty$ (or $\mathbb{K}$-*smooth*) if it is of class $C_{\mathbb{K}}^k$ for all $k \in \mathbb{N}_0$. If $\mathbb{K}$ is understood, we simply write $C^k$ instead of $C_{\mathbb{K}}^k$, and we call $f$ smooth or of class $C^\infty$ if it is $\mathbb{K}$-smooth.

**1.3** For example, every continuous linear mapping $\lambda \colon E \to F$ is smooth, with $\lambda^{[1]}(x, y, t) = \lambda(y)$ for all $(x, y, t) \in E \times E \times \mathbb{K}$. If $V, W$ and $F$ are topological $\mathbb{K}$-vector spaces and $\beta \colon V \times W \to F$ is a continuous bilinear map, then $\beta$ is smooth, with
$$\beta^{[1]}((v, w), (v', w'), t) = \beta(v, w') + \beta(v', w) + t\beta(v', w')$$
for all $v, v' \in V$, $w, w' \in W$, and $t \in \mathbb{K}$ (cf. [2]).

**1.4** If $k \geq 2$, then a map $f$ is of class $C_{\mathbb{K}}^k$ if and only if $f$ is of class $C_{\mathbb{K}}^1$ and $f^{[1]}$ is of class $C_{\mathbb{K}}^{k-1}$; in this case, $f^{[k]} = (f^{[1]})^{[k-1]}$ [2, Rem. 4.2].

**1.5** Given a $C^1_{\mathbb{K}}$-map $f\colon U \to F$ as before, we define
$$df(x,v) := \lim_{0 \neq t \to 0} \tfrac{1}{t}(f(x+tv) - f(x)) = f^{[1]}(x,v,0)$$
for $(x,v) \in U \times E$. Then $df\colon U \times E \to F$ is continuous, being a partial map of $f^{[1]}$, and it can be shown that the "differential" $df(x, \bullet)\colon E \to F$ of $f$ at $x$ is a continuous $\mathbb{K}$-linear map, for each $x \in U$ [2, Prop. 2.2]. If $f$ is of class $C^2$, we define a continuous mapping $d^2 f\colon U \times E^2 \to F$ via
$$d^2 f(x, v_1, v_2) := \lim_{0 \neq t \to 0} \tfrac{1}{t}(d(x + tv_2, v_1) - df(x, v_1)) = f^{[2]}((x, v_1, 0), (v_2, 0, 0), 0).$$
Similarly, if $f$ is of class $C^k_{\mathbb{K}}$, we obtain continuous mappings $d^j f \colon U \times E^j \to F$ for all $j \in \mathbb{N}_0$ such that $j \leq k$ (where $d^0 f := f$). It can be shown that $d^j f(x; \bullet)\colon E^j \to F$ is a symmetric $j$-multilinear map [2, La. 4.8].

Our discussion of implicit function theorems in the real case will be made easy by the following fact ([2, Prop. 7.4]):

**Proposition 1.6** *Let $E$ be a real topological vector space, $F$ a locally convex real topological vector space, $U \subseteq E$ an open subset, $f\colon U \to F$ a map, and $k \in \mathbb{N}_0 \cup \{\infty\}$. Then $f$ is of class $C^k_{\mathbb{R}}$ if and only if it is a $C^k$-map in the sense of Michal-Bastiani, i.e., $f$ is continuous, the differentials $d^j f \colon U \times E^j \to F$ described in **1.5** exist for all $j \in \mathbb{N}$ such that $j \leq k$,[1] and are continuous.* □

In more general situations, it is necessary to work with the functions $f^{[j]}$, since the differentials alone do not encode enough information. For example, $d^j f = 0$ for all $j \geq 2$ if $\mathbb{K}$ is a non-discrete topological field of characteristic 2 and $f$ any smooth function on $\mathbb{K}$ (cf. [2, Thm. 5.4]). Even worse, injective smooth functions $f\colon \mathbb{Z}_p \to \mathbb{Q}_p$ are known to exist whose derivative vanishes identically [32, Exercise 29.G].

**1.7** (Chain Rule). If $E$, $F$, and $H$ are topological $\mathbb{K}$-vector spaces, $U \subseteq E$ and $V \subseteq F$ are open subsets, and $f\colon U \to V \subseteq F$, $g\colon V \to H$ are mappings of class $C^k$, then also the composition $g \circ f \colon U \to H$ is of class $C^k$. If $k \geq 1$, we have $(f(x), f^{[1]}(x,y,t), t) \in V^{[1]}$ for all $(x,y,t) \in U^{[1]}$, and
$$(g \circ f)^{[1]}(x,y,t) = g^{[1]}(f(x), f^{[1]}(x,y,t), t). \tag{2}$$
In particular, $d(g \circ f)(x,y) = dg(f(x), df(x,y))$ for all $(x,y) \in U \times E$ [2, Prop. 3.1 and 4.5].

We recall from [2, La. 4.9] that being of class $C^k$ is a local property.

**Lemma 1.8** *Let $E$ and $F$ be topological $\mathbb{K}$-vector spaces, and $f \colon U \to F$ be a mapping, defined on an open subset $U$ of $E$. Let $k \in \mathbb{N}_0 \cup \{\infty\}$. If there is an open cover $(U_i)_{i \in I}$ of $U$ such that $f|_{U_i}\colon U_i \to F$ is of class $C^k$ for each $i \in I$, then $f$ is of class $C^k$.* □

---
[1] That is to say, all of the limits occurring in the recursive definition of the differentials $d^j f$ exist.

**Definition 1.9** A *valued field* is a field $\mathbb{K}$, together with an absolute value $|.|\colon \mathbb{K} \to [0, \infty[$ (see [38]); we require furthermore that the absolute value be non-trivial (meaning that it gives rise to a non-discrete topology on $\mathbb{K}$). An *ultrametric field* is a valued field $(\mathbb{K}, |.|)$ whose absolute value satisfies the ultrametric inequality

$$|x+y| \leq \max\{|x|, |y|\} \quad \text{for all } x, y \in \mathbb{K}.$$

Locally compact, totally disconnected, non-discrete topological fields will be referred to as *local fields*.

**Remark 1.10** It is well-known that every local field $\mathbb{K}$ admits an ultrametric absolute value defining its topology [37]. Fixing such an absolute value on $\mathbb{K}$, we can consider $\mathbb{K}$ as an ultrametric field.

**Remark 1.11** Note that we do not require that valued fields (nor ultrametric fields) be complete (with respect to the metric induced by the absolute value). Whenever our results depend on completeness properties of the ground field, we will state these explicitly.

**1.12** Recall that a topological vector space $E$ over an ultrametric field $\mathbb{K}$ is called *locally convex* if every zero-neighbourhood of $E$ contains an open $\mathbb{O}$-submodule of $E$, where $\mathbb{O} := \{t \in \mathbb{K} : |t| \leq 1\}$ is the valuation ring of $\mathbb{K}$. Equivalently, $E$ is locally convex if and only if its vector topology is defined by a family of ultrametric continuous seminorms $\gamma \colon E \to [0, \infty[$ on $E$ (cf. [29] for more information). Let $\mathbb{K}$ be a valued field. We call a topological $\mathbb{K}$-vector space *polynormed* if its vector topology is defined by a family of continuous seminorms (which need not be ultrametric when $\mathbb{K}$ is an ultrametric field). This terminology deviates from the one in Bourbaki [5], where only polynormed vector spaces over ultrametric fields are considered whose topology arises from a family of continuous *ultrametric* seminorms, and which therefore are precisely the locally convex spaces over such fields in our terminology. Ultrametric seminorms are called "ultra-semi-norms" in [5] and [6]. Occasionally, we shall write $\|.\|_\gamma$ for a continuous seminorm $\gamma$.

**1.13** A *Banach space* over a valued field $\mathbb{K}$ is a normed $\mathbb{K}$-vector space $(E, \|.\|)$ (see [6, Ch. I, §1, no. 2]) which is complete in the metric associated with $\|.\|$. Given a normed $\mathbb{K}$-vector space $(E, \|.\|)$, $x \in E$ and $r > 0$, we let $B_r^E(x) := \{y \in E \colon \|y - x\| < r\}$ be the open ball of radius $r$ around $x$. We write $B_r(x) := B_r^E(x)$ when no confusion is possible. $\overline{B}_r(x) := \{y \in E \colon \|y - x\| \leq r\}$ denotes the corresponding closed ball.

**1.14** We shall not presume that normed spaces (nor Banach spaces) over ultrametric fields be ultrametric, unless saying so explicitly. For example, $\ell^1(\mathbb{Q}_p)$ is a non-ultrametric (and non-locally convex) Banach space over $\mathbb{Q}_p$.

The following fact ([2, Thm. 12.4]) is essential for our discussion of implicit function theorems over ultrametric fields (cf. [21] for the real locally convex case):

**Proposition 1.15** *Let $(\mathbb{K}, |.|)$ be either $\mathbb{R}$, equipped with the usual absolute value, or an ultrametric field. Let $E$ and $F$ be topological $\mathbb{K}$-vector spaces (which need not be locally convex when $\mathbb{K} = \mathbb{R}$), $f: U \to F$ be a mapping, defined on a non-empty open subset $U \subseteq E$, and $k \in \mathbb{N}_0$. If $E$ is metrizable, then the following conditions are equivalent:*

(a) *$f$ is a mapping of class $C^k_\mathbb{K}$.*

(b) *The composition $f \circ c : \mathbb{K}^{k+1} \to F$ is of class $C^k_\mathbb{K}$, for every smooth mapping $c: \mathbb{K}^{k+1} \to U$.*

*In particular, $f$ is smooth if and only if $f \circ c$ is smooth, for every $k \in \mathbb{N}$ and every smooth map $c: \mathbb{K}^k \to U$.* $\square$

Finally, let us describe the notions of real and complex analytic mappings used in the present paper, and their basic properties.

We recall from [3], Definition 5.6:

**Definition 1.16** Let $E$ be a complex topological vector space, $F$ be a locally convex complex topological vector space, $U \subseteq E$ be an open subset, and $f: U \to F$ be a map. Then $f$ is called *complex analytic* if it is continuous and for every $x \in U$, there exists a zero-neighbourhood $V \subseteq E$ such that $x + V \subseteq U$ and

$$f(x+h) = \sum_{n=0}^{\infty} \beta_n(h) \quad \text{for all } h \in V$$

as a pointwise limit, for suitable continuous homogeneous polynomials $\beta_n : E \to F$ of degree $n \in \mathbb{N}_0$.

Real analytic mappings are defined as follows:

**Definition 1.17** Let $E$ be a real topological vector space, $F$ be a locally convex real topological vector space, $U \subseteq E$ be open, and $f: U \to F$ be a map. Then $f$ is called *real analytic* if it extends to a complex analytic mapping $V \to F_\mathbb{C}$, defined on some open neighbourhood $V$ of $U$ in $E_\mathbb{C}$.

**Remark 1.18** Real analyticity of a mapping $f: U \to F$ is a local property in the sense that real analyticity of $f|_{U_j}$ for an open cover $(U_j)_{j \in J}$ of $U$ entails real analyticity of $f$. In fact, if $f$ is real analytic locally, then for every $x \in U$ we find an open, balanced zero-neighbourhood $W_x \subseteq E$ such that $x + W_x \subseteq U$ and $f|_{x+W_x} = g_x|_{x+W_x}$ for some complex analytic mapping $g_x : V_x := (x + W_x) + iW_x \to F_\mathbb{C}$. Given $x, y \in U$, we have $V_{x,y} := V_x \cap V_y = ((x + W_x) \cap (y + W_y)) + i(W_x \cap W_y)$, where $W_x \cap W_y$ is a balanced open zero-neighbourhood and thus connected. Consequently, the connected components of $V_{x,y}$ are of the form $C + i(W_x \cap W_y)$, where $C$ is a connected component of $(x+W_x) \cap (y+W_y)$. Since $g_x$ and $g_y$ coincide on $C$, they coincide on $C + i(W_x \cap W_y)$ by the Identity Theorem [3, Prop. 6.6 II]. Thus $g_x|_{V_{x,y}} = g_y|_{V_{x,y}}$ for all $x, y \in U$, whence $g := \bigcup_{x \in U} g_x : \bigcup_{x \in U} V_x \to F_\mathbb{C}$ is a well-defined complex analytic mapping extending $f$.

As a consequence of [2, Prop. 7.7 and La. 10.1], we have:

**Lemma 1.19** *Every real or complex analytic map $f\colon U \to F$ as before is of class $C^\infty_{\mathbb{R}}$.* $\square$

For the next observation, see [2, Prop. 7.7]:

**Lemma 1.20** *Let $E$ be a complex topological vector space, $F$ be a locally convex complex topological vector space, $U \subseteq E$ be open, and $f\colon U \to F$ be a map. Then $f$ is complex analytic if and only if $f$ is of class $C^\infty_{\mathbb{C}}$, if and only if $f$ is of class $C^\infty_{\mathbb{R}}$ and $df(x,\bullet)\colon E \to F$ is complex linear for each $x \in U$.* $\square$

The Chain Rule for $C^\infty_{\mathbb{C}}$-functions readily entails that compositions of composable complex analytic (resp., real analytic) mappings are complex analytic (resp., real analytic).

## 2 Generalized Implicit Function Theorem in the Real and Complex Case

We begin with a preparatory result, providing *continuous* implicit functions in very general situations. A mapping from an open subset of a normed real vector space to a real locally convex space will be called an $FC^k$-map if it is $k$ times continuously differentiable in the Fréchet sense (cf. [5, 2.3.1]).

**Proposition 2.1** *Let $k \in \mathbb{N} \cup \{\infty\}$, $P$ be a topological space, $E$ a real Banach space, and $U \subseteq E$ an open subset. Suppose that $f\colon P \times U \to E$ is a continuous function such that*

(i) *$f_p := f(p,\bullet)\colon U \to E$ is of class $FC^k$, for all $p \in P$, and*

(ii) *The map $P \times U \to L(E)$, $(p,x) \mapsto f'_p(x) := d(f_p)(x,\bullet) = d_2 f(p,x,\bullet)$ is continuous, where $L(E)$ is equipped with the operator norm.*

*Let $(p,x) \in P \times U$, and suppose that $A := f'_p(x) \in \mathrm{GL}(E) := L(E)^\times$. Let $0 < a < 1 < b$ be given. Then there exists an open neighbourhood $Q \subseteq P$ of $p$ and $r > 0$ such that the open ball $B := B_r(x)$ of radius $r$ around $x$ is contained in $U$ and the following holds:*

(a) *$f_q(B)$ is open in $E$, for each $q \in Q$, and $\phi_q\colon B \to f_q(B)$, $\phi_q(y) := f_q(y) = f(q,y)$ is an invertible $FC^k$-map, with inverse $(\phi_q)^{-1}\colon f_q(B) \to B$ of class $FC^k$.*

(b) *For all $q \in Q$, $y \in B$, and $s \in \,]0, r - \|y - x\|]$, we have*

$$f_q(y) + A.B_{as}(0) \subseteq f_q(B_s(y)) \subseteq f_q(y) + A.B_{bs}(0). \tag{3}$$

(c) *$W := \bigcup_{q \in Q}(\{q\} \times f_q(B))$ is open in $P \times E$ and the map $\psi\colon W \to B$, $\psi(q,v) := \phi_q^{-1}(v)$ is continuous. Furthermore, the map $\theta\colon Q \times B \to W$, $\theta(q,y) := (q, f(q,y))$ is a homeomorphism, with inverse given by $\theta^{-1}(q,v) = (q, \psi(q,v))$.*

(d) *$Q \times (f_p(x) + A.B_\delta(0)) \subseteq W$ for some $\delta > 0$.*

*In particular, for each $q \in Q$, there is a unique element $\beta(q) \in B$ such that $f(q, \beta(q)) = f(p, x)$, and the mapping $\beta\colon Q \to B$ so obtained is continuous.*

**Proof.** In view of hypotheses (i) and (ii), we find an open neighbourhood $Q_1 \subseteq P$ of $p$ and $R > 0$ such that $B_R(x) \subseteq U$ and the following holds:

1. $f'_q(y) \in \mathrm{GL}(E)$ for all $q \in Q_1$ and $y \in B_R(x)$;
2. $\|f'_{q_1}(y_1)^{-1} f'_{q_2}(y_2) - \mathbf{1}\| < \frac{1}{2}$ for all $q_1, q_2 \in Q_1$ and $y_1, y_2 \in B_R(x)$;
3. $\|f'_q(y)^{-1}(f'_q(y_1) - f'_q(y_2))\| \leq 1 - \sqrt{a}$ for all $q \in Q_1$ and $y, y_1, y_2 \in B_R(x)$; and
4. $\|A^{-1} f'_q(y)\| < b$ and $\|f'_q(y)^{-1} A\| \leq \frac{1}{\sqrt{a}}$ for all $q \in Q_1$ and $y \in B_R(x)$.

Choose $r \in {]}0, \frac{R}{2}[$. There is an open neighbourhood $Q \subseteq Q_1$ of $p$ such that
$$\|A^{-1}.(f_q(x) - f_p(x))\| < \tfrac{ar}{2} \quad \text{for all } q \in Q. \tag{4}$$
Define $\delta := \frac{ar}{2}$. We claim that all assertions of the proposition are satisfied with $Q$, $r$, and $\delta$ as just defined.

(a) Given $q \in Q$, we consider the map $h\colon B \to E$, $h(y) := A^{-1} f_q(y)$. Then $\|h'(y) - \mathbf{1}\| = \|A^{-1} f'_q(y) - \mathbf{1}\| < \frac{1}{2}$ by 2., whence $h$ is injective. In fact, if $y_1 \neq y_2 \in B$, then
$$\begin{aligned}
\|h(y_2) - h(y_1)\| &\geq \|y_2 - y_1\| - \|h(y_2) - y_2 - (h(y_1) - y_1)\| \\
&= \|y_2 - y_1\| - \left\|\int_0^1 (h'(y_1 + t(y_2 - y_1)) - \mathbf{1}).(y_2 - y_1)\, dt\right\| \\
&\geq \tfrac{1}{2}\|y_2 - y_1\| > 0.
\end{aligned}$$

Hence also $f_q|_B = A \circ h$ is injective, and since $f'_q(y) \in \mathrm{GL}(E)$ for all $y \in B$ (by 1.) and $f_q$ is $FC^k$, the standard Inverse Function Theorem [23, Thm. I.5.2] shows that $f_q(B)$ is open in $E$ and $\phi_q^{-1}\colon f_q(B) \to B$ an $FC^k$-map, where $\phi_q := f_q|_B^{f_q(B)}$.

(b) Let $y \in B$, $s \in {]}0, r - \|y - x\|{]}$, and $q \in Q$. Given $z \in B_s(y)$, we have
$$\max\{\|A^{-1}.f'_q(y + t(z - y))\| : t \in [0, 1]\} < b$$
by 4., entailing that
$$\|A^{-1}.(f_q(z) - f_q(y))\| = \left\|\int_0^1 A^{-1}.f'_q(y + t(z - y)).(z - y)\, dt\right\| < bs.$$

Thus $f_q(z) \in f_q(y) + A.B_{bs}(0)$, verifying the second inclusion in (3).

To tackle the first inclusion in (3), let $y$, $s$, $q$ be as before; we want to apply [23, La. I.5.4] (with $\rho$, $s'$ as below playing the role of $r, s$) to the $FC^1$-map
$$g\colon B_s(0) \to E, \quad g(v) := f'_q(y)^{-1}(f_q(y + v) - f_q(y)).$$

Let $z \in B_{as}(0)$. Then $z' := f'_q(y)^{-1}.A.z \in B_{\frac{as}{\sqrt{a}}}(0) = B_{s\sqrt{a}}(0)$, by 4. If we can find $w \in B_s(0)$ such that $g(w) = z'$, then $y + w \in B_s(y)$ and

$$f_q(y+w) = f'_q(y).g(w) + f_q(y) = f_q(y) + A.z,$$

showing that $f_q(y) + A.z \in f_q(B_s(y))$, as desired. Now, since $\|z'\| < s\sqrt{a}$, we find $\rho \in ]\frac{\|z'\|}{\sqrt{a}}, s[$. Then $\|z'\| < \rho\sqrt{a} = (1-s')\rho$ with $s' := 1 - \sqrt{a}$. Since $\rho < s$, we have $\overline{B}_\rho(0) \subseteq B_s(0)$. Furthermore, $\|g'(v_2) - g'(v_1)\| = \|f'_q(y)^{-1}.(f'_q(y+v_2) - f'_q(y+v_1))\| \leq 1 - \sqrt{a} = s'$ for all $v_1, v_2 \in \overline{B}_\rho(0)$, by 3. Hence [23, La. I.5.4] provides $w \in \overline{B}_\rho(y) \subseteq B_s(0)$ such that $g(w) = z'$, as required. Thus, also the first inclusion in (3) holds.

(c) Let $q \in Q$, $y \in B$, and $\varepsilon \in ]0, r - \|y-x\|]$. There is an open neighbourhood $V \subseteq Q$ of $q$ such that $f(q,y) - f(q',y) \in A.B_{\frac{a\varepsilon}{2}}(0)$ for all $q' \in V$. Given $q' \in V$ and $z \in f_q(y) + A.B_{\frac{a\varepsilon}{2}}(0) \subseteq f_{q'}(y) + A.B_{a\varepsilon}(0)$, by (b) there exists $w \in B_\varepsilon(y) \subseteq B$ such that $\phi_{q'}(w) = f_{q'}(w) = z$, whence $(q', z) \in W$ and $w = \phi_{q'}^{-1}(z) \in B_\varepsilon(y)$. Thus $Y := V \times (f_q(y) + A.B_{\frac{a\varepsilon}{2}}(0))$ is an open neighbourhood of $(q, f_q(y))$, contained in $W$. We conclude that $W$ is open. Furthermore, for all $(q', z)$ in the open neighbourhood $Y$ of $(q, f_q(y))$ we have $\psi(q', z) = \phi_{q'}^{-1}(z) \in B_\varepsilon(y) = B_\varepsilon(\psi(q, f_q(y)))$. The continuity of $\psi$ follows. Now, apparently $\theta$ is continuous and is a bijection whose inverse has the asserted form. Hence, the map $\psi$ being continuous, so is $\theta^{-1}$.

(d) Let $q \in Q$. Applying (b) with $y := x$ and $s := r$, we get $f_q(x) + A.B_{ar}(0) \subseteq f_q(B)$. Therefore $f_p(x) + A.B_\delta(0) = f_p(x) + A.B_{\frac{ar}{2}}(0) = (f_p(x) - f_q(x)) + f_q(x) + A.B_{\frac{ar}{2}}(0) \subseteq A.B_{\frac{ar}{2}}(0) + f_q(x) + A.B_{\frac{ar}{2}}(0) = f_q(x) + A.B_{ar}(0) \subseteq f_q(B)$. Thus $\{q\} \times (f_p(x) + A.B_\delta(0)) \subseteq W$.

The final conclusion is clear; we have $\beta(q) = \psi(q, f(p,x)) = \phi_q^{-1}(f(p,x))$. □

**Lemma 2.2** *Let $P$ be an open subset of a real topological vector space $Z$, $E$ be a real Banach space, $U \subseteq E$ be open, $f: P \times U \to E$ be a map, and $k \in \mathbb{N} \cup \{\infty\}$. If $f$ is of class $C^{k+1}$ or if $E$ is finite-dimensional and $f$ is of class $C^k$, then hypotheses* (i) *and* (ii) *of Proposition* 2.1 *are satisfied. Furthermore, the mapping*

$$h: P \times U \to L(E), \quad h(p,x) := f'_p(x) = df((p,x), (0, \bullet))$$

*is of class $C^{k-1}$.*

**Proof.** If $E$ is finite-dimensional and $f$ is of class $C^k$, then $f_p := f(p, \bullet): U \to E$ is a $C^k$-map between open subsets of a finite-dimensional space, hence $k$ times continuously partially differentiable in the traditional sense, and thus an $FC^k$-map, as is well-known. Let $e_1, \ldots, e_n$ be a basis of $E$. The mappings $P \times U \to E$, $(p,x) \mapsto f'_p(x).e_j = d(f_p)(x, e_j) = df((p,x), (0, e_j))$ being of class $C^{k-1}$ for $j = 1, \ldots, n$, we readily deduce that $P \times U \to L(E) \cong M_n(\mathbb{R})$, $(p,x) \mapsto f'_p(x)$ is $C^{k-1}$ and hence continuous.

If $E$ is infinite-dimensional and $f$ is of class $C^{k+1}$, then, for each $p \in P$, $f_p := f(p, \bullet): U \to E$ is a $C^{k+1}$-map (in the Michal-Bastiani sense) between open subsets of Banach spaces and therefore an $FC^k$-map (see [19, Cor. 2.7.2 and p. 110], or, for a direct proof, [11, appendix]). The continuous linear map $\lambda\colon E \to Z \times E$, $\lambda(y) := (0, y)$ gives rise to a continuous linear (and hence smooth) map

$$L(\lambda, E)\colon L(Z \times E, E) \to L(E), \quad A \mapsto A \circ \lambda,$$

where $L(Z \times E, E)$ is equipped with the topology of uniform convergence on bounded sets. The mapping $h$ can be written as the composition $h = L(\lambda, E) \circ f'$, where

$$f'\colon P \times U \to L(Z \times E, E), \quad f'(p, x) := df((p, x), \bullet)$$

is of class $C^{k-1}$ by [11, Prop. 2.1] (which remains valid for non-locally convex domains, with identical proof). Hence $h$ is of class $C^{k-1}$ (and thus continuous). $\square$

**Theorem 2.3 (Generalized Implicit Function Theorem)** *Let $k \in \mathbb{N} \cup \{\infty\}$, $\mathbb{K} \in \{\mathbb{R}, \mathbb{C}\}$, $Z$ be a topological $\mathbb{K}$-vector space, $P \subseteq Z$ an open subset, $E$ a Banach space over $\mathbb{K}$, $U \subseteq E$ an open subset, and $f\colon P \times U \to E$ a map. We consider two situations:*

(i) *$\mathbb{K} = \mathbb{R}$ and $f$ is of class $C_\mathbb{R}^{k+1}$; or $\mathbb{K} = \mathbb{R}$, $E$ is finite-dimensional, and $f$ is of class $C_\mathbb{R}^k$; respectively:*

(ii) *$f$ is $\mathbb{K}$-analytic.*

*Let $(p, x) \in P \times U$, and suppose that $A := f_p'(x) \in \mathrm{GL}(E)$, where $f_p := f(p, \bullet)$. Furthermore, let $0 < a < 1 < b$ be given. Then there exists an open neighbourhood $Q \subseteq P$ of $p$ and $r > 0$ such that $B := B_r(x) \subseteq U$ and the following holds:*

(a) *$f_q(B)$ is open in $E$, for each $q \in Q$, and $\phi_q\colon B \to f_q(B)$, $\phi_q(y) := f_q(y) = f(q, y)$ is an invertible $FC^k$-map (resp., $\mathbb{K}$-analytic map), whose inverse $(\phi_q)^{-1}\colon f_q(B) \to B$ is of class $FC^k$ (resp., $\mathbb{K}$-analytic).*

(b) *For all $q \in Q$, $y \in B$, and $s \in {]}0, r - \|y - x\|]$, we have*

$$f_q(y) + A.B_{as}(0) \subseteq f_q(B_s(y)) \subseteq f_q(y) + A.B_{bs}(0).$$

(c) *$W := \bigcup_{q \in Q} (\{q\} \times f_q(B))$ is open in $Z \times E$, and the map $\psi\colon W \to B$, $\psi(q, v) := \phi_q^{-1}(v)$ is of class $C_\mathbb{R}^k$ (resp., $\mathbb{K}$-analytic). Furthermore, the map $\theta\colon Q \times B \to W$, $\theta(q, y) := (q, f(q, y))$ is a $C_\mathbb{R}^k$-diffeomorphism (resp., a $\mathbb{K}$-analytic diffeomorphism), with inverse given by $\theta^{-1}(q, v) = (q, \psi(q, v))$.*

(d) *$Q \times (f_p(x) + A.B_\delta(0)) \subseteq W$ for some $\delta > 0$.*

*In particular, for each $q \in Q$ there is a unique element $\beta(q) \in B$ such that $f(q, \beta(q)) = f(p, x)$, and the mapping $\beta\colon Q \to B$ so obtained is of class $C_\mathbb{R}^k$ (resp., $\mathbb{K}$-analytic).*

**Proof.** Let $f$, $p$, $x$, $a$, and $b$ be given as described in the theorem. Then hypotheses (i) and (ii) of Proposition 2.1 are satisfied, by Lemma 2.2. We let $Q$, $r$, $B$, $W$, $\psi$, $\theta$, and $\delta$ be as described in Proposition 2.1. Then (b) and (d) of the theorem hold by Proposition 2.1 (b) and (d), and in view of Part (a) of the proposition, apparently Part (a) of the theorem will hold if we can establish (c). We already know that $W$ is open, and we know that $\psi$ is continuous.

Let us assume first that we are in the situation of (i), and prove by induction on $j \in \mathbb{N}$, $j \leq k$ that $\psi$ is of class $C_\mathbb{R}^j$. If $j = 1$, suppose that $(q, v) \in W$ and $(q_1, v_1) \in Z \times E$. Set $y := \phi_q^{-1}(v)$. There is $\tau > 0$ such that $(q, v) + ]-\tau, \tau[\,(q_1, v_1) \subseteq W$. Define $g \colon ]-\tau, \tau[\, \times U \to \mathbb{R} \times E$ via $g(t, u) := (t, f(q + tq_1, u))$. Then $g$ is a mapping between open subsets of the Banach space $\mathbb{R} \times E$ which is of class $C_\mathbb{R}^k$ if $E$ (and thus $\mathbb{R} \times E$) is finite-dimensional, otherwise of class $C_\mathbb{R}^{k+1}$. In any case, $g$ is an $FC^k$-map and thus an $FC^1$-map, and apparently $g'(0, y)$ is invertible, as it can be considered as a lower triangular $2 \times 2$ block matrix with invertible diagonal entries $\mathrm{id}_\mathbb{R}$ and $f_q'(y)$. By the classical Inverse Function Theorem for Banach spaces ([23], Theorem I.5.2), we find $0 < \sigma \leq \min\{\tau, r - \|y - x\|\}$ such that $g$ restricts to an $FC^1$-diffeomorphism $h$ from $]-\sigma, \sigma[\, \times B_\sigma(y)$ onto an open subset $S \subseteq \mathbb{R} \times E$. There is $0 < \kappa < \sigma$ and an open neighbourhood $V \subseteq E$ of $v$ such that $]-\kappa, \kappa[\, \times V \subseteq S$. Then apparently $h^{-1}(t, w) = (t, \phi_{q+tq_1}^{-1}(w)) = (t, \psi(q + tq_1, w))$ for all $(t, w) \in ]-\kappa, \kappa[\, \times V$, entailing that $]-\kappa, \kappa[\, \times V \to E$, $(t, w) \mapsto \psi(q + tq_1, w)$ is an $FC^1$-map. After shrinking $\kappa$, we may assume that $v + ]-\kappa, \kappa[\, v_1 \subseteq V$. Then, by the preceding, $c \colon ]-\kappa, \kappa[\, \to E$, $t \mapsto \psi(q + tq_1, v + tv_1)$ is of class $C_\mathbb{R}^1$, and thus $d\psi((q, v), (q_1, v_1)) = c'(0)$ exists. Since $f(q + tq_1, c(t)) = v + tv_1$, the Chain Rule and the Rule on Partial Derivatives give

$$v_1 = \tfrac{d}{dt}|_{t=0} f(q + tq_1, c(t)) = d_1 f(q, \psi(q, v); q_1) + d_2 f(q, \psi(q, v); d\psi((q, v), (q_1, v_1))) \, .$$

Now $d_2 f(q, \psi(q, v); \bullet) = f_q'(\psi(q, v))$ is invertible by Proposition 2.1 (a). Thus

$$\begin{aligned} d\psi((q, v), (q_1, v_1)) &= f_q'(\psi(q, v))^{-1}.(v_1 - d_1 f(q, \psi(q, v); q_1)) \\ &= \varepsilon\big(f_q'(\psi(q, v))^{-1}, v_1 - d_1 f(q, \psi(q, v); q_1)\big) \end{aligned} \qquad (5)$$

for all $(q, v) \in W$ and $(q_1, v_1) \in Z \times E$, where $\varepsilon \colon L(E) \times E \to E$ is the bilinear evaluation map, which is a continuous since $E$ is normed. Now $\varepsilon$, $\psi$, $d_1 f$, inversion $\iota \colon \mathrm{GL}(E) \to \mathrm{GL}(E)$ and the mapping

$$P \times U \to L(E), \quad h(s, u) := f_s'(u) \qquad (6)$$

being continuous, we deduce from (5) that $d\psi \colon W \times Z \times E \to E$ is continuous, whence $\psi$ is of class $C_\mathbb{R}^1$. Similarly, if $1 < j < k$ and $\psi$ is of class $C_\mathbb{R}^j$ by induction hypothesis, using that $\varepsilon$ (being continuous bilinear), $\iota$ (cf. [12]) and the map in (6) (by Lemma 2.2) are of class $C_\mathbb{R}^j$, we deduce from (5) and the Chain Rule that $d\psi$ is of class $C_\mathbb{R}^j$. Thus $\psi$ is $C_\mathbb{R}^1$ with $d\psi$ of class $C_\mathbb{R}^j$, and hence $\psi$ is of class $C_\mathbb{R}^{j+1}$ (cf. [10]). Thus, the assertions of the theorem hold in situation (i).

We may pass to the complex analytic and real analytic cases using standard ideas from the finite-dimensional case (cf. [8, (10.2.4)]). Indeed, if $\mathbb{K} = \mathbb{C}$ and $f$ is complex analytic, then $f$ is of class $C_\mathbb{R}^\infty$ in particular, and thus $\psi$ is of class $C_\mathbb{R}^\infty$, by what has been shown. Equation (5) shows that $d\psi(w, \bullet) : Z \times E \to Z \times E$ is complex linear, for all $w \in W$. With Lemma 1.20, we deduce that $\psi$ is complex analytic. Finally, assume that $\mathbb{K} = \mathbb{R}$ and assume that $f$ is real analytic. Equip $E_{\mathbb{C}} \cong E \times E$ with the maximum norm. Given $(q, z) \in W$, set $y := \psi(q, z)$. There is a complex analytic function $F : Y \to E_{\mathbb{C}}$, defined on an open neighbourhood $Y$ of $P \times U$ in $Z_{\mathbb{C}} \times E_{\mathbb{C}}$, such that $F|_{P \times U} = f$. Then $d_2 F(q, y; \bullet) = f'_q(y)_{\mathbb{C}}$ being invertible, by the complex analytic case of the theorem just established, there exist open neighbourhoods $P_1 \subseteq Z_{\mathbb{C}}$ of $q$ and $r_1 > 0$ such that $P_1 \times B_1 \subseteq Y$ with $B_1 := B_{r_1}^{E_{\mathbb{C}}}(y)$, and such that $F(p_1, \bullet)|_{B_1}$ is a complex analytic diffeomorphism onto an open set, for each $p_1 \in P_1$, and furthermore $W_1 := \bigcup_{p_1 \in P_1} \{p_1\} \times F(\{p_1\} \times B_1)$ is open in $Z_{\mathbb{C}} \times E_{\mathbb{C}}$, and $\psi_1 : W_1 \to E_{\mathbb{C}}$, $\psi_1(p_1, v_1) := (F(p_1, \bullet)|_{B_1})^{-1}(v_1)$ is complex analytic. Then $Q \cap P_1$ and $B \cap B_1$ are open neighbourhoods of $q$ and $y$ in $Z$ and $E$, respectively. The map $\theta$ being a homeomorphism onto the open set $W$, we deduce that $W_2 := \bigcup_{s \in Q \cap P_1} (\{s\} \times f_s(B \cap B_1)) = \theta((Q \cap P_1) \times (B \cap B_1))$ is an open subset of $Z \times E$. Then $\psi_1 : W_1 \to E_{\mathbb{C}}$ is a complex analytic mapping on an open neighbourhood of $W_2$ in $Z_{\mathbb{C}} \times E_{\mathbb{C}}$ which extends $\psi|_{W_2}$, and thus $\psi|_{W_2}$ is real analytic. Being real analytic locally by the preceding, $\psi$ is real analytic (Remark 1.18). $\square$

## 3 Strict Differentiability

We now leave the framework of real and complex analysis and turn to differential calculus over arbitrary valued fields. In order to be able to prove implicit function theorems, we require a differentiability property which is stronger than being $C^1$, namely "strict differentiability." In this section, we recall the definition of strictly differentiable mappings from open subsets of normed $\mathbb{K}$-vector spaces to polynormed $\mathbb{K}$-vector spaces, where $\mathbb{K}$ is a valued field. We show that every strictly differentiable mapping is of class $C^1$, and we show that, conversely, every mapping of class $C^2$ from an open subset of a normed space to a polynormed $\mathbb{K}$-vector space is strictly differentiable.

**Definition 3.1** Let $\mathbb{K}$ be a valued field, $E$ be a normed $\mathbb{K}$-vector space, $F$ a polynormed $\mathbb{K}$-vector space, $U \subseteq E$ be open, and $f : U \to F$ be a map. Given $x \in U$, we say that $f$ is *strictly differentiable at $x$* if there exists a continuous linear map $A \in L(E, F)$ such that, for every $\varepsilon > 0$ and continuous seminorm $\gamma$ on $F$, there exists $\delta > 0$ such that

$$\|f(z) - f(y) - A.(z-y)\|_\gamma < \varepsilon \, \|z - y\|$$

for all $y, z \in U$ such that $\|z - x\| < \delta$ and $\|y - x\| < \delta$. The map $f$ is called *strictly differentiable* if it is strictly differentiable at each $x \in U$.

It is clear that $A$ is uniquely determined in the preceding situation; we write $f'(x) := A$.

If $E$ is a normed vector space over a valued field $\mathbb{K}$, and $F$ a polynormed $\mathbb{K}$-vector space, we equip the space $L(E, F)$ of continuous $\mathbb{K}$-linear maps $E \to F$ with the topology of uniform convergence on bounded subsets of $E$. This topology makes $L(E, F)$ a polynormed $\mathbb{K}$-vector space, whose vector topology arises from the family of continuous seminorms

$$\|A\|_\gamma := \sup\{\|A.v\|_\gamma \cdot \|v\|^{-1} \colon 0 \neq v \in E\} \in [0, \infty[$$

(cf. [30], p. 59), where $\gamma$ ranges through the continuous seminorms on $F$. If also $F$ is normed, with norm $\gamma$, then $L(E, F)$ is normable; its vector topology arises from the operator norm $\|.\| := \|.\|_\gamma$.

**Lemma 3.2** *Let $\mathbb{K}$ be a valued field, $E$ be a normed $\mathbb{K}$-vector space, $F$ a polynormed $\mathbb{K}$-vector space, $U \subseteq E$ be open, and $f \colon U \to F$ be a strictly differentiable map. Then $f$ is of class $C^1$, we have $f'(x) = df(x, \bullet)$ for all $x \in U$, and the mapping*

$$f' \colon U \to L(E, F), \quad x \mapsto f'(x)$$

*is continuous.*

**Proof.** *Directional derivatives.* Given $x \in U$ and $y \in E$, let us show that the directional derivative $df(x, y)$ exists, and is given by $f'(x).y$. For $y = 0$ this is trivial. If $0 \neq y \in E$, there exists $r > 0$ such that $x + ty \in U$ for all $0 \neq t \in B_r(0) \subseteq \mathbb{K}$. By strict differentiability of $f$ in $x$, for every continuous seminorm $\gamma$ on $F$ we have

$$\left\|\tfrac{1}{t}(f(x+ty) - f(x)) - f'(x)y\right\|_\gamma = \|y\| \cdot \frac{\|f(x+ty) - f(x) - f'(x)ty\|_\gamma}{\|ty\|} \to 0$$

as $t \to 0$, showing that $df(x, y) := \lim_{0 \neq t \to 0} \tfrac{1}{t}(f(x+ty) - f(x)) = f'(x).y$ indeed.

*$f'$ is continuous.* In fact, given $\varepsilon > 0$ and a continuous seminorm $\gamma$ on $F$, we find $\delta > 0$ such that $\|f(z) - f(y) - f'(x).(z - y)\|_\gamma < \varepsilon \|z - y\|$ for all $y, z \in U$ such that $\|y - x\| < \delta$ and $\|z - x\| < \delta$. Let $y \in U$ such that $\|y - x\| < \delta$. Then, given $0 \neq u \in E$ we have $y + tu \in U$ and $\|y + tu - x\| < \delta$ for all $t \in \mathbb{K}^\times$ sufficiently close to $0$, entailing that

$$\begin{aligned}\|(f'(y) - f'(x)).u\|_\gamma &= \lim_{t \to 0} \left\|\tfrac{1}{t}(f(y+tu) - f(y)) - f'(x).u\right\|_\gamma \\ &= \|u\| \cdot \lim_{t \to 0} \frac{\|f(y+tu) - f(y) - f'(x).(tu)\|_\gamma}{\|tu\|} \leq \|u\| \cdot \varepsilon\,.\end{aligned}$$

As a consequence, $\|f'(y) - f'(x)\|_\gamma \leq \varepsilon$ for all $y \in U$ such that $\|y - x\| < \delta$, showing that $f' \colon U \to L(E, F)$ is continuous.

*$f$ is of class $C^1_\mathbb{K}$.* Note first that $f$ is continuous. In fact, given $x \in U$, $\varepsilon > 0$, and a continuous seminorm $\gamma$ on $F$, we find $\delta > 0$ as in Definition 3.1. Pick $0 < \rho \leq \delta$ such that $\rho \cdot (\varepsilon + \|f'(x)\|_\gamma) < \varepsilon$. Then, for all $y \in B_\rho(x) \cap U$, we estimate

$$\begin{aligned}\|f(y) - f(x)\|_\gamma &\leq \|f(y) - f(x) - f'(x).(y - x)\|_\gamma + \|f'(x).(y - x)\|_\gamma \\ &\leq \varepsilon \cdot \|y - x\| + \|f'(x)\|_\gamma \cdot \|y - x\| < \varepsilon\,.\end{aligned}$$

We deduce that $f$ is continuous.

Next, let $W := \{(x,y,t) \in U^{[1]} : t \neq 0\}$. Define $g: U^{[1]} \to F$ via $g(x,y,t) := \frac{1}{t}(f(x+ty) - f(x))$ for $(x,y,t) \in W$, $g(x,y,0) := f'(x).y$ for $(x,y) \in U \times E$. Then $g|_W$ is continuous since $f$ is continuous, and $g|_{U \times E \times \{0\}}: U \times E \times \{0\} \to F$ is continuous since $f'$ and the evaluation map $L(E,F) \times E \to F$ are continuous (the space $E$ being normed). Hence $g$ will be continuous if we can show that $g(x_\alpha, y_\alpha, t_\alpha) \to g(x,y,0)$, for every net $((x_\alpha, y_\alpha, t_\alpha))_{\alpha \in I}$ in $W$ which converges to some $(x,y,0) \in U^{[1]}$. Since $\|y_\alpha\| \leq \|y\| + 1$ eventually, we may assume without loss of generality that $\|y_\alpha\| \leq \|y\| + 1$ for all $\alpha$. If $y \neq 0$, then $y_\alpha \neq 0$ eventually, whence we may furthermore assume in this case that $y_\alpha \neq 0$ for all $\alpha$. Then, for a given $\alpha$, we either have $y_\alpha = 0$ (in which case $y = 0$ by the preceding): then
$$\|g(x_\alpha, y_\alpha, t_\alpha) - g(x,y,0)\| = \|0 - 0\| = 0\,;$$
or we have $y_\alpha \neq 0$, in which case

$$\begin{aligned}
&\|g(x_\alpha, y_\alpha, t_\alpha) - g(x,y,0)\|_\gamma \\
&= \left\|\tfrac{1}{t_\alpha}(f(x_\alpha - t_\alpha y_\alpha) - f(x_\alpha)) - f'(x).y\right\|_\gamma \\
&\leq \|y_\alpha\| \cdot \frac{\|f(x_\alpha + t_\alpha y_\alpha) - f(x_\alpha) - f'(x).t_\alpha y_\alpha\|_\gamma}{|t_\alpha| \cdot \|y_\alpha\|} + \|f'(x).(y_\alpha - y)\|_\gamma \\
&\leq (\|y\| + 1) \cdot \frac{\|f(x_\alpha + t_\alpha y_\alpha) - f(x_\alpha) - f'(x).t_\alpha y_\alpha\|_\gamma}{\|t_\alpha y_\alpha\|} + \|f'(x)\|_\gamma \cdot \|y_\alpha - y\|\,,
\end{aligned}$$

where the first term tends to 0 as $\alpha$ increases since $f$ is strictly differentiable at $x$, and the second term tends to 0 for trivial reasons. $\square$

We want to show that every $C^2$-map is strictly differentiable. The proof hinges on symmetry properties of $f^{[1]}$ and $f^{[2]}$, as described in following lemma:

**Lemma 3.3** *Suppose that $E$ and $F$ are topological vector spaces over a non-discrete topological field $\mathbb{K}$, and $f: U \to F$ is a map, defined on an open subset of $E$.*

(a) *If $f$ is $C^1$, $t \in \mathbb{K}^\times$, and $(x,y,s) \in E \times E \times \mathbb{K}$ such that $(x,y,ts) \in U^{[1]}$, then also $(x,ty,s) \in U^{[1]}$, and*
$$t\,f^{[1]}(x,y,ts) = f^{[1]}(x,ty,s)\,. \tag{7}$$

(b) *If $f$ is $C^2$, $t \in \mathbb{K}^\times$, $x, x_1, y, y_1 \in E$ and $s, s_1, s_2 \in \mathbb{K}$ such that*
$$((x,y,ts), (x_1, y_1, ts_1), ts_2) \in U^{[2]}\,,$$
*then also $((x, t^2 y, \tfrac{s}{t}), (tx_1, t^3 y_1, s_1), s_2) \in U^{[2]}$, and*
$$t^3\,f^{[2]}((x,y,ts), (x_1, y_1, ts_1), ts_2) = f^{[2]}((x, t^2 y, \tfrac{s}{t}), (tx_1, t^3 y_1, s_1), s_2)\,. \tag{8}$$

**Proof.** (a) Since $x + (ts)y = x + s(ty)$, it is obvious that $(x, ty, s) \in U^{[1]}$ if and only if $(x, y, ts) \in U^{[1]}$. In this case, we have
$$tf^{[1]}(x, y, ts) = \tfrac{1}{s}(f(x + tsy) - f(x)) = f^{[1]}(x, ty, s)$$
provided $s \neq 0$; if $s = 0$, then $f^{[1]}(x, ty, s) = f^{[1]}(x, ty, 0) = df(x, ty) = tdf(x, y) = tf^{[1]}(x, y, 0) = tf^{[1]}(x, y, ts)$. Thus (7) is established.

(b) Let $t \in \mathbb{K}^\times$, $x, y, x_1, y_1 \in E$, and $s, s_1, s_2 \in \mathbb{K}$ such that $((x, y, ts), (x_1, y_1, ts_1), ts_2) \in U^{[2]}$. Then $(x, y, ts) \in U^{[1]}$ and hence also $(x, t^2y, \tfrac{s}{t}) \in U^{[1]}$, by Part (a). If $s_2 = 0$, this entails that $((x, t^2y, \tfrac{s}{t}), (tx_1, t^3 y_1, s_1), s_2) \in U^{[2]}$. If $s_2 \neq 0$, we calculate:

$$\begin{aligned}
f^{[2]}((x, y, ts), (x_1, y_1, ts_1), ts_2) &= (f^{[1]})^{[1]}((x, y, ts), (x_1, y_1, ts_1), ts_2) \\
&= \tfrac{1}{ts_2}(f^{[1]}(x + ts_2 x_1, y + ts_2 y_1, \underbrace{ts + t^2 s_1 s_2}_{=t^2(\tfrac{s}{t}+s_1 s_2)}) - f^{[1]}(x, y, \underbrace{ts}_{=t^2 \cdot \tfrac{s}{t}})) \\
&= \tfrac{1}{t^3 s_2}\left(f^{[1]}(x + ts_2 x_1, t^2 y + t^3 s_2 y_1, \tfrac{s}{t} + s_1 s_2) - f^{[1]}(x, t^2 y, \tfrac{s}{t})\right) \\
&= \tfrac{1}{t^3} f^{[2]}((x, t^2 y, \tfrac{s}{t}), (tx_1, t^3 y_1, s_1), s_2),
\end{aligned}$$

showing that $((x, t^2 y, \tfrac{s}{t}), (tx_1, t^3 y_1, s_1), s_2) \in U^{[2]}$ and that (8) holds, when $s_2 \neq 0$. Here, we used Part (a) to pass to the third line. Letting $s_2 \neq 0$ tend to $0$, in view of the continuity of the functions involved we see that (8) remains valid for $s_2 = 0$. $\square$

Cf. also [12, La. 6.6] for a (less explicit) result concerning $f^{[k]}$ for arbitrary $k$.

We are now in the position to prove:

**Proposition 3.4** *Let $\mathbb{K}$ be a valued field, $E$ be a normed $\mathbb{K}$-vector space, $F$ a polynormed $\mathbb{K}$-vector space, $U \subseteq E$ be open, and $f: U \to F$ be a mapping of class $C^2$. Then $f$ is strictly differentiable.*

**Proof.** Given $x \in U$, let us show that $f$ is strictly differentiable at $x$. For all $y, z \in U$, we have

$$\begin{aligned}
&f(z) - f(y) - f'(x).(z - y) \\
&\quad = f(z) - f(y) - f'(y).(z - y) + f'(y).(z - y) - f'(x).(z - y) \\
&\quad = f^{[1]}(y, z - y, 1) - f^{[1]}(y, z - y, 0) + f^{[1]}(y, z - y, 0) - f^{[1]}(x, z - y, 0) \\
&\quad = f^{[2]}((y, z - y, 0), (0, 0, 1), 1) + f^{[2]}((x, z - y, 0), (y - x, 0, 0), 1). \tag{9}
\end{aligned}$$

Let us have a closer look at the individual terms. For each $t \in \mathbb{K}^\times$, we have

$$\begin{aligned}
f^{[2]}((y, z - y, 0), (0, 0, 1), 1) &= f^{[2]}((y, t^2 \cdot \tfrac{1}{t^2}(z - y), \tfrac{1}{t} \cdot 0), (t \cdot 0, t^3 \cdot 0, 1), 1) \\
&= t^3 f^{[2]}((y, \tfrac{1}{t^2}(z - y), 0), (0, 0, t), t), \tag{10}
\end{aligned}$$

using Lemma 3.3 (b). The function $f^{[2]}((x, \bullet, 0), (y - x, 0, 0), 1) = f'(y) - f'(x) \colon E \to F$ is linear. Furthermore, we have

$$f^{[2]}((x, z - y, 0), (y - x, 0, 0), 1) = s f^{[2]}((x, z - y, 0), (\tfrac{1}{s}(y - x), 0, 0), s) \tag{11}$$

for all $s \in \mathbb{K}^\times$, by Lemma 3.3 (a). Let $\varepsilon > 0$ and $\gamma$ be a continuous seminorm on $F$. Since $f^{[2]}((x,0,0), (0,0,0), 0) = d^2f(x,0,0) = 0$, in view of the continuity of $f^{[2]}$ and openness of $U^{[2]}$, there exists $\delta > 0$ such that $((u,v,0), (w,0,a), b) \in U^{[2]}$ and

$$\|f^{[2]}((u,v,0), (w,0,a), b)\|_\gamma < \varepsilon \quad \text{for all } u \in B_\delta^E(x),\ v, w \in B_\delta^E(0),\ a, b \in B_\delta^{\mathbb{K}}(0). \quad (12)$$

We may assume that $\delta \leq 1$. Pick $\rho \in \mathbb{K}^\times$ such that $|\rho| < 1$, and $s \in \mathbb{K}^\times$ such that $|s| \leq \frac{\delta |\rho|^2}{2}$. Define $r := \min\{|s|\delta, \frac{1}{8}\delta^3|\rho|^6\}$. Let $y, z \in B_r^E(x) \subseteq U$ such that $y \neq z$. Then there is a unique integer $k \in \mathbb{Z}$ such that

$$|\rho|^{k+1} \leq \sqrt{\frac{\|z - y\|}{\delta}} < |\rho|^k.$$

Set $t := \rho^k$. Then $\|t^{-2}(z-y)\| = |\rho|^{-2k}\|z-y\| < \delta$,

$$|t| = |\rho|^k \leq \frac{1}{|\rho|} \cdot \sqrt{\frac{\|z-y\|}{\delta}} \leq \frac{1}{|\rho|} \cdot \sqrt{\frac{2r}{\delta}} \leq \frac{1}{|\rho|}\sqrt{\frac{1}{4}\frac{\delta^3|\rho|^6}{\delta}} = \frac{\delta|\rho|^2}{2} < \delta,$$

and

$$\frac{|t|^3}{\|z-y\|} = |\rho|^{k-2}\frac{|\rho|^{2k+2}}{\|z-y\|} \leq \frac{|\rho|^{k-2}}{\delta} \leq \frac{1}{\delta|\rho|^3}\sqrt{\frac{\|z-y\|}{\delta}} \leq \frac{1}{\delta|\rho|^3}\sqrt{\frac{2r}{\delta}} \leq \frac{1}{\delta|\rho|^3}\sqrt{\frac{1}{4}\delta^2|\rho|^6} = \frac{1}{2}.$$

Hence, using (10) and (12) (with $u := y$),

$$\frac{\|f^{[2]}((y, z-y, 0), (0,0,1), 1)\|_\gamma}{\|z-y\|} = \frac{|t|^3}{\|z-y\|}\|f^{[2]}((y, \tfrac{1}{t^2}(z-y), 0), (0,0,t), t)\|_\gamma < \tfrac{\varepsilon}{2}. \quad (13)$$

Using (11) with $s$ as just chosen, we obtain

$$\frac{\|f^{[2]}((x, z-y, 0), (y-x, 0, 0), 1)\|_\gamma}{\|z-y\|}$$
$$= \frac{|s|}{\|z-y\|}\|f^{[2]}((x, z-y, 0), (\tfrac{1}{s}(y-x), 0, 0), s)\|_\gamma$$
$$= \frac{|s| \cdot |\rho|^{2k}}{\|z-y\|}\|f^{[2]}((x, \rho^{-2k}(z-y), 0), (\tfrac{1}{s}(y-x), 0, 0), s)\|_\gamma < \tfrac{\varepsilon}{2}. \quad (14)$$

Indeed, we have $|s| < \delta$ by choice of $s$, $|s|^{-1}\|y-x\| < |s|^{-1}r \leq \delta$, and $|\rho|^{-2k}\|z-y\| < \delta$, whence $\|f^{[2]}((x, \rho^{-2k}(z-y), 0), (\tfrac{1}{s}(y-x), 0, 0), s)\|_\gamma < \varepsilon$ holds, by (12). Furthermore, $|s||\rho|^{2k}\|z-y\|^{-1} \leq \frac{|s|}{\delta|\rho|^2} \leq \frac{1}{2}$. By (9), (13) and (14), we have

$$\frac{\|f(z) - f(y) - f'(x).(z-y)\|_\gamma}{\|z-y\|} < \varepsilon$$

for all $z \neq y \in B_r(x) \subseteq U$. Thus $f$ is strictly differentiable at $x$. $\square$

The following variant of Proposition 3.4 involving parameters will be needed later:

**Lemma 3.5** *Let $(\mathbb{K}, |.|)$ be a valued field, $E$ be a normed $\mathbb{K}$-vector space, $U \subseteq E$ an open subset, $F$ be a polynormed $\mathbb{K}$-vector space, $Z$ a topological $\mathbb{K}$-vector space, and $P \subseteq Z$ an open subset. Let $f\colon P \times U \to F$ be a mapping of class $C^2_{\mathbb{K}}$. Then the map*

$$h\colon P \times U \to L(E, F), \quad (p, x) \mapsto f'_p(x) \tag{15}$$

*is continuous, where $f_p := f(p, \bullet)$ and $f'_p(x) := d(f_p)(x, \bullet)$. For all $p \in P$, $x \in U$, $\varepsilon > 0$ and continuous seminorm $\gamma$ on $F$, there exists a neighbourhood $Q$ of $p$ in $P$ and $r > 0$ such that*

$$\frac{\|f_q(z) - f_q(y) - f'_q(x).(z-y)\|_\gamma}{\|z - y\|} < \varepsilon \quad \text{for all } q \in Q \text{ and } y \neq z \in B_r(x) \cap U. \tag{16}$$

**Proof.** Let $p \in P$, $x \in U$, $\varepsilon > 0$ and a continuous seminorm $\gamma$ on $F$ be given. Since

$$g\colon P \times U^{[2]} \to F, \quad (q, z) \mapsto (f_q)^{[2]}(z)$$

is continuous as a partial map of $f^{[2]}$ and $g(p, (x, 0, 0), (0, 0, 0), 0) = d^2(f_p)(x, 0, 0) = 0$, there exists $\delta > 0$ and an open neighbourhood $Q$ of $p$ in $P$ such that $((u, v, 0), (w, 0, a), b) \in U^{[2]}$ and

$$\|(f_q)^{[2]}((u, v, 0), (w, 0, a), b)\|_\gamma = \|g(q, (u, v, 0), (w, 0, a), b)\|_\gamma < \varepsilon,$$

for all $q \in Q$, $u \in B^E_\delta(x)$, $v, w \in B^E_\delta(0)$, and $a, b \in B^{\mathbb{K}}_\delta(0)$. For each fixed $q \in Q$, replacing $f$ with $f_q$ in the preceding proof, we find $r$ (independent of $q$) as described there and can repeat the estimates verbatim, to obtain (16).

To see that $h$ is continuous, let $p \in P$, $x \in U$, $\varepsilon > 0$ and a continuous seminorm $\gamma$ on $F$ be given. Since $f^{[2]}(((p, x), (0, 0), 0), ((0, 0), (0, 0), 0), 0) = d^2 f((p, x), (0, 0), (0, 0)) = 0$ and $f^{[2]}$ is continuous, there is $\delta > 0$ and a balanced zero-neighbourhood $V \subseteq Z$ such that

$$\|f^{[2]}(((p, x), (0, u), 0), ((v, z), (0, 0), 0), t)\|_\gamma \leq \varepsilon$$

for all $u, z \in B^E_\delta(0)$, $v \in V$, and $t \in B^{\mathbb{K}}_\delta(0)$. Pick $\rho \in \mathbb{K}^\times$ such that $|\rho| < 1$, and pick $t \in \mathbb{K}^\times$ such that $|t| \leq \delta|\rho|$. Define $r := \delta|t|$. We claim that

$$\|f'_q(y) - f'_p(x)\|_\gamma \leq \varepsilon \tag{17}$$

for all $q \in P \cap (p + tV)$ and $y \in B_r(x) \cap U$. To see this, let $0 \neq u \in E$. There is a unique $k \in \mathbb{Z}$ such that $|\rho|^{k+1} \leq \frac{\|u\|}{\delta} < |\rho|^k$. Then, using Lemma 3.3 (a) and linearity in $u$,

$$\begin{aligned}
\frac{\|(f'_q(y) - f'_p(x)).u\|_\gamma}{\|u\|} &= \frac{\|f^{[2]}(((p, x), (0, u), 0), ((q-p, y-x), (0, 0), 0), 1)\|_\gamma}{\|u\|} \\
&= \tfrac{|t|\,|\rho|^k}{\|u\|} \|f^{[2]}(((p, x), (0, \rho^{-k}u), 0), ((\tfrac{1}{t}(q-p), \tfrac{1}{t}(y-x)), (0, 0), 0), t)\|_\gamma \\
&\leq \varepsilon\,,
\end{aligned}$$

entailing that (17) holds. Hence $h$ is continuous, which completes the proof. $\square$

# 4 Uniform Differentiability

For mappings on open subsets of finite-dimensional topological vector spaces over locally compact topological fields with values in polynormed spaces, the results of the preceding section can be strengthened substantially: such a mapping is of class $C^1$ if and only if it is strictly differentiable, if and only if it is "locally uniformly differentiable."

**Definition 4.1** Suppose that $(\mathbb{K}, |.|)$ is a valued field, $(E, \|.\|)$ a normed $\mathbb{K}$-vector space, $F$ a polynormed $\mathbb{K}$-vector space, $U \subseteq E$ an open subset, and $f: U \to F$ a map. Then $f$ is called *uniformly differentiable* if there exists a function $f': U \to L(E, F)$ such that, for every $\varepsilon > 0$ and continuous seminorm $\gamma$ on $F$, there exists $\delta > 0$ with the following property: for all $x, y, z \in U$ such that $\|y - x\| < \delta$, $\|z - x\| < \delta$ and $y \neq z$, we have

$$\frac{\|f(z) - f(y) - f'(x).(z - y)\|_\gamma}{\|z - y\|} < \varepsilon \,.$$

We call $f$ *locally uniformly differentiable* if every $x \in U$ has an open neighbourhood $V \subseteq U$ such that $f|_V$ is uniformly differentiable.

**Remark 4.2** It is apparent from the definitions that every locally uniformly differentiable mapping is strictly differentiable.

**Remark 4.3** Strengthening Lemma 3.2, it is easy to see that $f': U \to L(E, F)$ is uniformly continuous, for every uniformly differentiable mapping $f: E \supseteq U \to F$. As we shall not need this fact, the simple proof is omitted.

**Lemma 4.4** *Let $\mathbb{K}$ be a locally compact, non-discrete topological field, and $|.|$ be an absolute value on $\mathbb{K}$ defining its topology. Let $E$ be a finite-dimensional normed $\mathbb{K}$-vector space, $F$ be a polynormed $\mathbb{K}$-vector space, $U \subseteq E$ be open, and $f: U \to F$ be a mapping of class $C^1$. Then $f$ is locally uniformly differentiable. If, furthermore, $U$ is compact, then $f$ is uniformly differentiable.*

**Proof.** (cf. [32], Exercise 28E (ii) when $E = F = \mathbb{K}$). Pick $0 \neq \rho \in \mathbb{K}$ such that $|\rho| < 1$. Let $V \subseteq U$ be an open subset with compact closure $\overline{V} \subseteq U$. Define $f': U \to L(E, F)$, $f'(x) := df(x, \bullet) = f^{[1]}(x, \bullet, 0)$. Given $\varepsilon > 0$ and a continuous seminorm $\gamma$ on $F$, consider the continuous function

$$g: U^{[1]} \to F, \quad g(x, y, t) := f^{[1]}(x, y, t) - f^{[1]}(x, y, 0) \,.$$

Then $\overline{V} \times \overline{B^E_{\frac{1}{|\rho|}}(0)} \times \{0\} \subseteq U \times E \times \{0\} \subseteq U^{[1]}$ is a compact subset on which $g$ vanishes identically. Using a compactness argument, we find $\sigma > 0$ such that $\overline{V} \times \overline{B^E_{\frac{1}{|\rho|}}(0)} \times B^{\mathbb{K}}_\sigma(0) \subseteq U^{[1]}$ and such that $\|g(x, y, t)\|_\gamma < \frac{\varepsilon}{2}$ for all $(x, y, t) \in \overline{V} \times \overline{B^E_{\frac{1}{|\rho|}}(0)} \times B^{\mathbb{K}}_\sigma(0)$. Let $e_1, \ldots, e_n$ be a basis of $E$, and $e_1^*, \ldots, e_n^* \in E'$ be its dual basis.

Note that, given $A \in L(E,F)$, for every $0 \neq v \in E$ we have $\frac{\|A.v\|_\gamma}{\|v\|} = \frac{\|\sum_{i=1}^n e_i^*(v)A.e_i\|_\gamma}{\|v\|} \leq \sum_{i=1}^n \frac{|e_i^*(v)|}{\|v\|} \cdot \|A.e_i\|_\gamma \leq \sum_{i=1}^n \|e_i^*\| \cdot \|A.e_i\|_\gamma$. Thus

$$\|A\|_\gamma \leq \sum_{i=1}^n \|e_i^*\| \cdot \|A.e_i\|_\gamma \text{ for all } A \in L(E,F). \tag{18}$$

Let $i \in \{1,\ldots,n\}$. The mapping $\overline{V} \to F$, $x \mapsto df(x,e_i)$ being continuous and thus uniformly continuous, we find $\delta_i > 0$ such that $\|df(y,e_i) - df(x,e_i)\|_\gamma < \frac{\varepsilon}{2n\|e_i^*\|}$ for all $x,y \in \overline{V}$ such that $\|x-y\| < \delta_i$. Define $\delta := \frac{1}{2}\min\{\sigma, \delta_1, \ldots, \delta_n\}$. By (18) and the choice of $\delta_i$, we have $\|df(y,\bullet) - df(x,\bullet)\|_\gamma < \frac{\varepsilon}{2}$ for all $x,y \in \overline{V}$ such that $\|x-y\| < 2\delta$.

Let $x,y,z \in V$ be given such that $y \neq z$, $\|y-x\| < \delta$, and $\|z-x\| < \delta$. There exists $k \in \mathbb{Z}$ such that $|\rho|^{k+1} \leq \|z-y\| < |\rho|^k$. We set $s := \rho^{k+1}$. Then $\|\frac{1}{s}(z-y)\| < \frac{1}{|\rho|}$, $|s| = |\rho|^{k+1} \leq \|z-y\| < 2\delta \leq \sigma$, and $\|z-y\| < 2\delta$. As a consequence,

$$\begin{aligned}
&\frac{\|f(z) - f(y) - f'(x).(z-y)\|_\gamma}{\|z-y\|} \\
&\leq \frac{\|f(z) - f(y) - f'(y).(z-y)\|_\gamma}{\|z-y\|} + \frac{\|(f'(y)-f'(x)).(z-y)\|_\gamma}{\|z-y\|} \\
&= \frac{|s|}{\|z-y\|} \cdot \left\|\tfrac{1}{s}(f(z)-f(y)) - f'(y).\tfrac{1}{s}(z-y)\right\|_\gamma + \tfrac{\varepsilon}{2} \\
&\leq \left\|f^{[1]}\left(y, \tfrac{1}{s}(z-y), s\right) - f^{[1]}\left(y, \tfrac{1}{s}(z-y), 0\right)\right\|_\gamma + \tfrac{\varepsilon}{2} \\
&= \left\|g\left(y, \tfrac{1}{s}(z-y), s\right)\right\|_\gamma + \tfrac{\varepsilon}{2} \leq \varepsilon.
\end{aligned}$$

We have shown that $f|_V$ is uniformly differentiable. The first assertion readily follows. To obtain the second assertion, choose $V := U$ in the first part of the proof. $\square$

We also need a variant of Lemma 3.5 involving parameters.

**Lemma 4.5** *Let $\mathbb{K}$ be a locally compact, non-discrete topological field, and $|.|$ be an absolute value on $\mathbb{K}$ defining its topology. Let $E$ be a finite-dimensional normed $\mathbb{K}$-vector space, $U \subseteq E$ be open, $F$ be a polynormed $\mathbb{K}$-vector space, and $P$ be a topological space. Let $f: P \times U \to F$ be a continuous mapping such that $f_p := f(p,\bullet): U \to F$ is of class $C^1_\mathbb{K}$ for all $p \in P$, and such that the mapping*

$$P \times U^{[1]} \to F, \quad (p,y) \mapsto (f_p)^{[1]}(y)$$

*is continuous. Let $p \in P$ and $u \in U$ be given. Then, for every $\varepsilon > 0$ and continuous seminorm $\gamma$ on $F$, there is a neighbourhood $Q$ of $p$ in $P$ and $\delta > 0$ such that*

$$\frac{\|f_q(z) - f_q(y) - f'_q(u).(z-y)\|_\gamma}{\|z-y\|} < \varepsilon$$

*for all $q \in Q$ and $y,z \in B_\delta(u) \cap U$ such that $y \neq z$, where $f'_q(u) := d(f_q)(u,\bullet)$.*

**Proof.** Let $\varepsilon > 0$ and $\gamma$ be given. Pick $0 \neq \rho \in \mathbb{K}$ such that $|\rho| < 1$. Let $V \subseteq U$ be an open neighbourhood of $u$ with compact closure $\overline{V} \subseteq U$. Consider the continuous mapping

$$g \colon P \times U^{[1]} \to F, \quad g(q, x, y, t) := f_q^{[1]}(x, y, t) - f_q^{[1]}(x, y, 0).$$

Then $\{p\} \times \overline{V} \times \overline{B^E_{\frac{1}{|\rho|}}(0)} \times \{0\} \subseteq P \times U \times E \times \{0\} \subseteq P \times U^{[1]}$ is a compact subset on which $g$ vanishes identically. Using a compactness argument, we find $\sigma > 0$ and a neighbourhood $P_0$ of $p$ in $P$ such that $\overline{V} \times \overline{B^E_{\frac{1}{|\rho|}}(0)} \times B^\mathbb{K}_\sigma(0) \subseteq U^{[1]}$ and such that

$$\|g(q, x, y, t)\|_\gamma < \tfrac{\varepsilon}{2} \quad \text{for all } (q, x, y, t) \in P_0 \times \overline{V} \times \overline{B^E_{\frac{1}{|\rho|}}(0)} \times B^\mathbb{K}_\sigma(0).$$

Let $e_1, \ldots, e_n$ be a basis of $E$, and $e_1^*, \ldots, e_n^*$ be its dual basis. Using the compactness of $\overline{V}$, we find a neighbourhood $Q \subseteq P_0$ of $p$ and $\kappa > 0$ such that $\|df_q(z, e_i) - df_q(y, e_i)\|_\gamma < \frac{\varepsilon}{2n \|e_i^*\|}$ for all $q \in Q$, $i \in \{1, \ldots, n\}$, and all $y, z \in \overline{V}$ such that $\|z - y\| < \kappa$. Define $\delta := \min\{\frac{\sigma}{2}, \frac{\kappa}{2}\}$. Re-using the estimates from the proof of Lemma 4.4, we see that the assertion of the lemma is satisfied for $Q$ and $\delta$. $\square$

## 5   Strict Differentiability of Higher Order

Generalizing the standard notion of (once) strictly differentiable mappings, in this section we define and discuss $k$ times strictly differentiable mappings on open subsets of normed vector spaces over valued fields.

**Definition 5.1** Let $\mathbb{K}$ be a valued field, $E$ be a normed $\mathbb{K}$-vector space, $F$ be a polynormed $\mathbb{K}$-vector space, and $U \subseteq E$ be an open subset. A function $f \colon U \to F$ is called an $SC^0$-*map* if it is continuous; it is called an $SC^1$-*map* is it is strictly differentiable (and hence $C^1$ in particular). Inductively, having defined $SC^k$-map for some $k \in \mathbb{N}$ (which are $C^k$ in particular), we call $f$ an $SC^{k+1}$-*map* if it is an $SC^k$-map and the mapping $f^{[k]} \colon U^{[k]} \to F$ is $SC^1$, where $E^{[k]}$ is equipped with the maximum norm. The map $f$ is called $SC^\infty$ if it is an $SC^k$-map for all $k \in \mathbb{N}_0$.

**Remark 5.2** In other words, $f$ is $SC^k$ if and only if $f$ is $C^k$ and $f^{[j]} \colon U^{[j]} \to F$ is strictly differentiable for all $j \in \mathbb{N}_0$ such that $j < k$. It follows from this and **1.4** that $f$ is $SC^k$ if and only if $f$ is $SC^1$ and $f^{[1]}$ is $SC^{k-1}$.

**Remark 5.3** If $f \colon E \supseteq U \to F$ of class $C^{k+1}$ in the preceding situation, then $f$ is an $SC^k$-map. In fact, for every $j \in \mathbb{N}_0$ such that $j < k$, the map $f^{[j]}$ is of class $C^{k+1-j}$, where $k + 1 - j \geq 2$, and hence strictly differentiable by Proposition 3.4.

**Remark 5.4** A mapping from an open subset of a finite-dimensional $\mathbb{K}$-vector space to a polynormed $\mathbb{K}$-vector space over a locally compact, non-discrete topological field $\mathbb{K}$ is of class $C^k$ if and only if it is an $SC^k$-map, by a simple induction based on Lemma 3.2, Lemma 4.4, and Remark 5.2.

Compositions of composable $SC^k$-maps are $SC^k$.

**Proposition 5.5** *Let $\mathbb{K}$ be a valued field, $E$ and $F$ be normed $\mathbb{K}$-vector spaces, $G$ be a polynormed $\mathbb{K}$-vector space, and $U \subseteq E$, $V \subseteq F$ be open subsets. Let $k \in \mathbb{N}_0 \cup \{\infty\}$ and suppose that $f\colon U \to V \subseteq F$ and $g\colon V \to G$ are $SC^k$-maps. Then $g \circ f\colon U \to G$ is an $SC^k$-map.*

**Proof.** We may assume that $k < \infty$; the proof is by induction. The case $k = 0$ is trivial. The case $k = 1$ is known (cf. [5, 1.3.1], where however no proof is given), and can be shown as follows: Given $x \in U$, let $\gamma$ be a continuous seminorm on $G$, and $\varepsilon > 0$. Set $\varepsilon' := \min\left\{\frac{1}{2(\|f'(x)\|+1)}, \frac{1}{2(\|g'(f(x))\|_\gamma+1)}\right\} < 1$. By strict differentiability of $g$ at $f(x)$, there exists $r > 0$ such that
$$\|g(z) - g(y) - g'(f(x)).(z-y)\|_\gamma \leq \varepsilon' \|z-y\| \quad \text{for all } y, z \in B_r^F(f(x)) \cap V.$$
By strict differentiability of $f$ at $x$ and continuity of $f$ at $x$, there exists $\delta > 0$ such that
$$\|f(z) - f(y) - f'(x).(z-y)\| \leq \varepsilon' \|z-y\| \quad \text{for all } z, y \in B_\delta^E(x) \cap U,$$
and $f(y) \in B_r^F(f(x))$ for all $y \in B_\delta^E(x) \cap U$. Then
$$\begin{aligned}
&\|g(f(z)) - g(f(y)) - g'(f(x)).f'(x).(z-y)\|_\gamma \\
&\leq \|g(f(z)) - g(f(y)) - g'(f(x)).(f(z)-f(y))\|_\gamma \\
&\quad + \|g'(f(x)).(f(z) - f(y) - f'(x).(z-y))\|_\gamma \\
&\leq \varepsilon' \underbrace{\|f(z) - f(y)\|}_{\leq (\|f'(x)\|+\varepsilon')\cdot\|z-y\|} + \varepsilon' \|g'(f(x))\|_\gamma \cdot \|z-y\| \\
&\leq \varepsilon \|z-y\|
\end{aligned}$$
for all $y, z \in B_\delta^E(x) \cap U$, whence $g \circ f$ is strictly differentiable at $x$, with differential $(g \circ f)'(x) = g'(f(x)) \circ f'(x)$.

*Induction step.* Assume that $2 \leq k \in \mathbb{N}$, and suppose that the assertion is correct when $k$ is replaced with $k-1$. Let $f$ and $g$ be $SC^k$-maps, as above. By the preceding, $g \circ f$ is an $SC^1$-map. We also know that
$$(g \circ f)^{[1]}(x,y,t) = g^{[1]}(f(x), f^{[1]}(x,y,t), t) \quad \text{for all } (x,y,t) \in U^{[1]}, \tag{19}$$
where $f^{[1]}$ and $g^{[1]}$ are $SC^{k-1}$-maps (cf. (2)). Now $f$ being an $SC^k$-map and hence an $SC^{k-1}$-map, and the continuous linear map $E \times E \times \mathbb{K} \to E$, $(x,y,t) \mapsto x$ being an $SC^{k-1}$-map, by induction the composition $U^{[1]} \to F$, $(x,y,t) \mapsto f(x)$ is an $SC^{k-1}$-map. Also $U^{[1]} \to \mathbb{K}$, $(x,y,t) \mapsto t$ is an $FC^{k-1}$-map, being the restriction of a continuous linear map. As a consequence,
$$\widehat{T}(f)\colon U^{[1]} \to V^{[1]}, \quad (x,y,t) \mapsto (f(x), f^{[1]}(x,y,t), t)$$
is an $SC^{k-1}$-map, its coordinates being $SC^{k-1}$ (cf. [2], proof of La. 4.4). But thus $(g \circ f)^{[1]} = g^{[1]} \circ \widehat{T}(f)$ is an $SC^{k-1}$-map, by the induction hypotheses. Now $g \circ f$ being $SC^1$ with $(g \circ f)^{[1]}$ being $SC^{k-1}$, we deduce that $g \circ f$ is an $SC^k$-map (see Remark 5.2). □

# 6 Newton Approximation with Parameters

In this section, we discuss Newton approximation and Newton approximation with parameters, as the basis for our inverse function theorems (resp., implicit function theorems) for valued fields.

**Lemma 6.1 (Newton approximation)** *Let $(\mathbb{K}, |.|)$ be a valued field, and $(E, \|.\|)$ be a Banach space over $\mathbb{K}$. Let $r > 0$, $x \in E$, and $f \colon B_r(x) \to E$ be a mapping. We suppose that there exists $A \in \mathrm{GL}(E) := L(E)^\times$ such that*

$$\sigma := \sup\left\{ \frac{\|f(z) - f(y) - A.(z-y)\|}{\|z - y\|} : y, z \in B_r(x),\ y \neq z \right\} < \frac{1}{\|A^{-1}\|}. \tag{20}$$

*Then the following holds:*

(a) *Let $a := 1 - \sigma \|A^{-1}\| \in\ ]0, 1]$ and $b := 1 + \sigma \|A^{-1}\| \in [1, 2[$. Then*

$$a\|z - y\| \leq \|A^{-1}.f(z) - A^{-1}.f(y)\| \leq b\|z - y\| \quad \text{for all } y, z \in B_r(x). \tag{21}$$

*For every $y \in B_r(x)$ and $s \in\ ]0, r - \|y\|]$, we have*

$$f(y) + A.B_{as}(0) \subseteq f(B_s(y)) \subseteq f(y) + A.B_{bs}(0). \tag{22}$$

*In particular, the map $f$ has open image and is a homeomorphism onto its image.*

(b) *If $(\mathbb{K}, |.|)$ is an ultrametric field in particular and $(E, \|.\|)$ an ultrametric Banach space, then $A^{-1} \circ f \colon B_r(x) \to E$ is isometric. For each $y \in B_r(x)$ and $s \in\ ]0, r]$, we have $B_s(y) \subseteq B_r(x)$ and*

$$f(B_s(y)) = f(y) + A.B_s(0). \tag{23}$$

**Proof.** (cf. [32, Lemma 27.4] when $E = \mathbb{K}$ is a complete ultrametric field. Compare also [20, Theorem 6.3.6] for a related result in the real case).

(a) Given $y, z \in B_r(x)$, we have

$$\begin{aligned}
\|A^{-1}.f(z) - A^{-1}.f(y)\| &= \|A^{-1}.(f(z) - f(y) - A.(z - y)) + z - y\| \\
&\leq \|A^{-1}\| \cdot \|f(z) - f(y) - A.(z - y)\| + \|z - y\| \\
&\leq (\sigma\|A^{-1}\| + 1)\|z - y\| = b\|z - y\|
\end{aligned}$$

and

$$\begin{aligned}
\|z - y\| &= \|A^{-1}.(f(z) - f(y) - A.(z - y)) - (A^{-1}.f(z) - A^{-1}.f(y))\| \\
&\leq \|A^{-1}\| \cdot \|f(z) - f(y) - A.(z - y)\| + \|A^{-1}.f(z) - A^{-1}.f(z)\| \\
&\leq \sigma\|A^{-1}\| \cdot \|z - y\| + \|A^{-1}.f(z) - A^{-1}.f(y)\|,
\end{aligned}$$

whence (21) holds. As a consequence of (21), $A^{-1} \circ f$ and hence also $f$ is injective and a homeomorphism onto its image. Now suppose that $y \in B_r(x)$ and $s \in {]0, r - \|y\|]}$. By the preceding, we have $f(B_s(y)) \subseteq f(y) + A.B_{bs}(0)$. To see that $f(y) + A.B_{as}(0) \subseteq f(B_s(y))$, let $c \in f(y) + A.B_{as}(0)$. There exists $t \in {]0,1[}$ such that $c \in f(y) + A.\overline{B}_{tas}(0)$. Given $z \in \overline{B}_{st}(y)$, we define
$$g(z) := z - A^{-1}.(f(z) - c).$$
Then $g(z) \in \overline{B}_{st}(y)$, since
$$\|g(z) - y\| \leq \underbrace{\|z - y - A^{-1}.f(z) + A^{-1}.f(y)\|}_{\leq \|A^{-1}\|\sigma\|z-y\| \leq \|A^{-1}\|\sigma st} + \underbrace{\|A^{-1}.c - A^{-1}.f(y)\|}_{\leq ats}$$
$$\leq (\|A^{-1}\|\sigma + a)st = st.$$
Thus $g(\overline{B}_{st}(y)) \subseteq \overline{B}_{st}(y)$. The map $g \colon \overline{B}_{st}(y) \to \overline{B}_{st}(y)$ is a contraction, since
$$\begin{aligned} \|g(u) - g(v)\| &= \|u - v - A^{-1}.(f(u) - f(v))\| \\ &\leq \|A^{-1}\| \cdot \|f(u) - f(v) - A.(u-v)\| \\ &\leq \sigma \cdot \|A^{-1}\| \cdot \|u - v\| \end{aligned} \qquad (24)$$
for all $u, v \in \overline{B}_{st}(y)$, where $\|A^{-1}\|\sigma < 1$. By Banach's Contraction Theorem ([32], p. 269), there exists a unique element $z_0 \in \overline{B}_{st}(y)$ such that $g(z_0) = z_0$ and hence $f(z_0) = z_0$. This completes the proof of (a).

(b) Now assume that $(\mathbb{K}, |.|)$ is an ultrametric field and $(E, \|.\|)$ an ultrametric Banach space. For all $y, z \in B_r(x)$ such that $y \neq z$, we have $\|A^{-1}.f(z) - A^{-1}.f(y) - (z-y)\| \leq \|A^{-1}\| \cdot \|f(z) - f(y) - A.(z-y)\| < \|z - y\|$, where we used (20) to obtain the final inequality. Hence, the norm $\|.\|$ being ultrametric, we must have $\|A^{-1}.f(z) - A^{-1}.f(y)\| = \|z - y\|$. Thus $A^{-1} \circ f$ is in fact isometric.

Given $y \in B_r(x)$ and $s \in {]0, r]}$, we have $B_s(y) \subseteq B_r(y) = B_r(x)$, exploiting that $\|.\|$ is ultrametric. Since $A^{-1} \circ f$ is isometric, we have $f(B_s(y)) = A.(A^{-1} \circ f)(B_s(y)) \subseteq A.B_s(A^{-1}.f(y)) = f(y) + A.B_s(0)$. If $c \in f(y) + A.B_s(0)$ is given, define
$$g(z) := z - A^{-1}.(f(z) - c) \text{ for } z \in B_s(y).$$
Then
$$\begin{aligned} \|g(z) - y\| &= \|(z-y) - (A^{-1}.f(z) - A^{-1}.f(y)) + A^{-1}.(c - f(y))\| \\ &\leq \max\{\|z-y\|, \|A^{-1}.f(z) - A^{-1}.f(y)\|, \|A^{-1}.(c - f(y))\|\} < s \end{aligned}$$
for $z \in B_s(y)$, whence $g(z) \in B_s(y)$. The map $g \colon B_s(y) \to B_s(y)$ is a contraction, by the calculation from (24). Recall that, the norm on $E$ being ultrametric, the open ball $B_s(y)$ is also closed and therefore complete in the induced metric. By Banach's Contraction Theorem ([32], p. 269), there is a unique element $z_0 \in B_s(y)$ such that $g(z_0) = z_0$ and thus $f(z_0) = c$. The proof is complete. $\square$

More generally, we shall need Newton approximation with parameters.

**Lemma 6.2 (Newton approximation with parameters)** *Let $(\mathbb{K}, |.|)$ be a valued field, $(E, \|.\|)$ be a Banach space over $\mathbb{K}$, and $P$ be a topological space. Let $r > 0$, $x \in E$, and $f\colon P \times B \to E$ be a continuous mapping, where $B := B_r^E(x)$. Given $p \in P$, we abbreviate $f_p := f(p, \bullet)\colon B \to E$. We suppose that there exists $A \in \mathrm{GL}(E) := L(E)^\times$ such that*

$$\sigma := \sup\left\{\frac{\|f_p(z) - f_p(y) - A.(z - y)\|}{\|z - y\|} : p \in P, \, y, z \in B, \, y \neq z\right\} < \frac{1}{\|A^{-1}\|}. \quad (25)$$

*Then $f_p(B)$ is open in $E$ and $f_p|_B$ is a homeomorphism onto its image, for each $p \in P$. The set $W := \bigcup_{p \in P}\{p\} \times f_p(B)$ is open in $P \times E$, and $\psi\colon W \to E$, $\psi(p, z) := (f_p|_B^{f_p(B)})^{-1}(z)$ is continuous. Furthermore, the map $\theta\colon P \times B \to W$, $\theta(p, y) := (p, f(p, y))$ is a homeomorphism, with inverse given by $\theta^{-1}(p, z) = (p, \psi(p, z))$.*

**Proof.** By Lemma 6.1, applied to $f_p$, the set $f_p(B)$ is open in $E$ and $f_p|_B$ a homeomorphism onto its image. Define $a := 1 - \sigma\|A^{-1}\|$. Let us show openness of $W$ and continuity of $h$. If $(p, z) \in W$, there exists $y \in B$ such that $f_p(y) = z$. Let $\varepsilon \in {]}0, r - \|y\|]$ be given. There is an open neighbourhood $Q$ of $p$ in $P$ such that $f(q, y) \in f(p, y) + A.B_{\frac{a\varepsilon}{2}}(0)$ for all $q \in Q$, by continuity of $f$. Then, as a consequence of Lemma 6.1 (a), Eqn. (22),

$$f_q(B_\varepsilon(y)) \supseteq f(q, y) + A.B_{a\varepsilon}(0) \supseteq f(p, y) + A.B_{\frac{a\varepsilon}{2}}(0) = z + A.B_{\frac{a\varepsilon}{2}}(0).$$

By the preceding, $Q \times (z + A.B_{\frac{a\varepsilon}{2}}(0)) \subseteq W$, whence $W$ is a neighbourhood of $(p, z)$. Furthermore, $\psi(q, z') = (f_q)^{-1}(z') \in B_\varepsilon(y) = B_\varepsilon((f_p)^{-1}(z)) = B_\varepsilon(\psi(p, z))$ for all $(q, z')$ in the neighbourhood $Q \times (z + A.B_{\frac{a\varepsilon}{2}}(0))$ of $(p, z)$. Thus $W$ is open and $\psi$ is continuous. The assertions concerning $\theta$ follow immediately. $\square$

# 7 Inverse Function and Implicit Function Theorems for $SC^k$-maps between Banach spaces

In this section, we prove an inverse function theorem and an implicit function theorem for $SC^k$-maps over complete valued fields, which parallel the classical theorems for continuously Fréchet differentiable mappings in the real case.

We begin with an Inverse Function Theorem for mappings strictly differentiable at a point (cf. also [5, 1.5.1]), strictly differentiable maps, and locally uniformly differentiable maps.

**Proposition 7.1** *Let $(\mathbb{K}, |.|)$ be a valued field, $(E, \|.\|)$ a Banach space over $\mathbb{K}$, $U \subseteq E$ an open subset, $x \in U$, and $f\colon U \to E$ be a mapping which is strictly differentiable at $x$ (resp., strictly differentiable, resp., locally uniformly differentiable), with $A := f'(x) \in \mathrm{GL}(E)$. Let $a, b \in \mathbb{R}$ be given such that $0 < a < 1 < b$. Then there exists $r > 0$ such that $B := B_r(x) \subseteq U$ and*

(a) $a\|z - y\| \leq \|A^{-1}.f(z) - A^{-1}.f(y)\| \leq b\|z - y\|$ *for all $y, z \in B$, whence $f|_B$ is injective in particular;*

(b) $f(y) + A.B_{as}(0) \subseteq f(B_s(y)) \subseteq f(y) + A.B_{bs}(0)$ for all $y \in B$ and $s \in ]0, r - \|y\|]$, whence, in particular, $f(B)$ is an open subset of $E$;

(c) The mapping $g := (f|_B^{f(B)})^{-1} : f(B) \to B$ is strictly differentiable at $f(x)$ (resp., strictly differentiable, resp., uniformly differentiable).

If $(\mathbb{K}, |.|)$ is an ultrametric field here and $(E, \|.\|)$ an ultrametric Banach space, then $r$ can be chosen such that (c) remains valid but (a) and (b) may be replaced with:

(a)' The mapping $A^{-1} \circ f|_B$ is isometric;

(b)' For all $y \in B$ and $s \in ]0, r]$, we have $f(B_s(y)) = f(y) + A.B_s(0)$. Thus, if $A$ is an isometry, then in fact $f(B_s(y)) = B_s(f(y))$.

**Proof.** Let $A := f'(x)$ and $c := \min\left\{\frac{b-1}{\|A^{-1}\|}, \frac{1-a}{\|A^{-1}\|}\right\}$. Then $1 - c\|A^{-1}\| \geq a$, $1 + c\|A^{-1}\| \leq b$, and $c < \frac{1}{\|A^{-1}\|}$. The map $f$ being strictly differentiable at $x$, there exists $r > 0$ such that $B_r(x) \subseteq U$ and

$$\frac{\|f(z) - f(y) - A.(z - y)\|}{\|z - y\|} \leq c \quad \text{for all } z, y \in B_r(x) \text{ such that } z \neq y.$$

Thus hypothesis (20) of Lemma 6.1 is satisfied; Parts (a), (b), (a)' and (b)' directly follow from that lemma.

(c) Assume that $f$ is locally uniformly differentiable, or strictly differentiable. Since $\mathrm{GL}(E)$ is open in $L(E)$, $f'(x) \in \mathrm{GL}(E)$, and $f'|_B$ is continuous, after shrinking $r$ we may assume that $f'(B_r(x)) \subseteq \mathrm{GL}(E)$. Inversion in $\mathrm{GL}(E)$ being continuous, after shrinking $r$ further we may also assume that $\|f'(y)^{-1}\| \leq \|A^{-1}\| + 1$ for all $y \in B_r(x)$.

*Assume that $f$ is locally uniformly differentiable.* After shrinking $r$ further if necessary, we may assume that $f$ is uniformly differentiable on $B := B_r(x)$. Abbreviate $g := (f|_B^{f(B)})^{-1}$ and define $g' : f(B) \to L(E)$ via $g'(y) := f'(g(y))^{-1}$. Given $\varepsilon > 0$, due to uniform differentiability of $f|_B$ there exists $\delta \in ]0, r]$ such that

$$\frac{\|f(w) - f(v) - f'(u).(w - v)\|}{\|w - v\|} < \frac{a\,\varepsilon}{(\|A^{-1}\| + 1)\|A^{-1}\|} \tag{26}$$

for all $u, v, w \in B$ such that $\|v - u\| < \delta$, $\|w - u\| < \delta$, and $v \neq w$. Set $\delta' := \frac{a\delta}{\|A^{-1}\|}$. Let $u', v', w' \in f(B)$ such that $v' \neq w'$, $\|v' - u'\| < \delta'$, and $\|w' - u'\| < \delta'$. Then $u := g(u')$, $v := g(v')$, and $w := g(w')$ are elements of $B$ such that $v \neq w$,

$$\begin{aligned}\|v - u\| &\leq a^{-1}\|A^{-1}.f(v) - A^{-1}.f(u)\| = a^{-1}\|A^{-1}.(v' - u')\| \\ &\leq a^{-1}\|A^{-1}\| \cdot \|v' - u'\| < a^{-1}\|A^{-1}\|\delta' = \delta\end{aligned}$$

(using Part (a)), and similarly $\|w - u\| < \delta$ and $\|w - v\| \le a^{-1}\|A^{-1}\| \cdot \|w' - v'\|$. Using (26) and Part (a), we obtain the following estimates:

$$\frac{\|g(w') - g(v') - g'(u').(w' - v')\|}{\|w' - v'\|} = \frac{\|w - v\|}{\|w' - v'\|} \cdot \frac{\|w - v - f'(u)^{-1}.(f(w) - f(v))\|}{\|w - v\|}$$

$$\le \|f'(u)^{-1}\| \cdot \frac{\|w - v\|}{\|w' - v'\|} \cdot \frac{\|f(w) - f(v) - f'(u).(w - v)\|}{\|w - v\|}$$

$$< (\|A^{-1}\| + 1) \cdot a^{-1} \cdot \|A^{-1}\| \cdot \frac{a\,\varepsilon}{(\|A^{-1}\| + 1)\|A^{-1}\|} = \varepsilon \,.$$

Thus $g$ is uniformly differentiable.

If $f$ is strictly differentiable at $x$, we see along the preceding lines (holding however $u := x$ and $u' := f(x)$ fixed now), that $g := (f|_B^{f(B)})^{-1}$ is strictly differentiable at $f(x)$, where $B := B_r(x)$.

If $f$ is strictly differentiable on all of $B$, because $f'(y) \in \mathrm{GL}(E)$ for all $y \in B$, we may apply the preceding proof just as well when $x$ is replaced with $y$; thus $g$ is strictly differentiable at each $z = f(y) \in f(B)$. □

Isometries are encountered frequently in the ultrametric case. If $\mathbb{K}$ is a valued field and $(E, \|.\|)$ a normed $\mathbb{K}$-vector space, we let $\mathrm{Iso}(E, \|.\|)$ denote the group of all bijective linear isometries of $E$. If the norm on $E$ is understood, we simply write $\mathrm{Iso}(E)$. We have:

**Lemma 7.2** *If $(E, \|.\|)$ is an ultrametric Banach space over an ultrametric field $\mathbb{K}$, then $\mathrm{Iso}(E)$ is open in $\mathrm{GL}(E) = L(E)^\times$.*

**Proof.** If $A \in L(E)$ such that $\|A\| < 1$, then $\|Ax\| < \|x\|$ and hence $\|(\mathbf{1} + A)x\| = \|x\|$ for all $0 \ne x \in E$, whence $\mathbf{1} + A$ is an isometry. Furthermore, using Neumann's series we see that $\mathbf{1} + A$ is invertible, with inverse $(\mathbf{1} + A)^{-1} = \sum_{k=0}^{\infty}(-1)^k A^k$. Thus $\mathbf{1} + A \in \mathrm{Iso}(E)$ for all $A \in L(E)$ such that $\|A\| < 1$, entailing that $\mathrm{Iso}(E)$ is open in $\mathrm{GL}(E)$. □

Let $\mathbb{K}$ be a valued field. An $SC^k$-*diffeomorphism* is an invertible $SC^k$-map $f \colon U \to V$ between open subsets of normed $\mathbb{K}$-vector spaces, such that $f^{-1}$ is an $SC^k$-map.

**Theorem 7.3 (Inverse Function Theorem for $SC^k$-maps)** *Let $k \in \mathbb{N} \cup \{\infty\}$, $(\mathbb{K}, |.|)$ be a valued field, $E$ be a Banach space over $\mathbb{K}$, $U \subseteq E$ be an open subset, and $f \colon U \to E$ be an $SC^k$-map. If $k > 1$, we assume that $\mathbb{K}$ is complete. Suppose that $f'(x) := df(x, \bullet) \in \mathrm{GL}(E)$ for some $x \in U$. Then there exists $r > 0$ such that $B := B_r(x) \subseteq U$, the set $f(B)$ is open in $E$, and $f|_B^{f(B)}$ is an $SC^k$-diffeomorphism.*

For our inductive proof of the Inverse Function Theorem, we need the Implicit Function Theorem, which we formulate as an "Inverse Function Theorem with Parameters."

**Theorem 7.4 (Implicit Function Theorem for $SC^k$-maps)** *Let $(\mathbb{K}, |.|)$ be a valued field, $k \in \mathbb{N} \cup \{\infty\}$, $Z$ and $F$ be Banach spaces over $\mathbb{K}$, $P \subseteq Z$ and $U \subseteq F$ be open subsets,* $f: P \times U \to F$ be an $SC^k$-map, and $(p, x) \in P \times U$ be a point such that $A := f'_p = d_2 f(p, x, \cdot) = df((p, x), (0, \cdot)) \in \mathrm{GL}(F)$; here $f_p := f(p, \cdot): U \to F$, and $E := Z \times F$ is equipped with the maximum norm, $\|(q, y)\| := \max\{\|q\|, \|y\|\}$ for all $q \in Z$, $y \in F$. If $k > 1$, we assume that $\mathbb{K}$ is complete. Let $a, b \in \mathbb{R}$ be given such that $0 < a < 1 < b$. Then there exists an open neighbourhood $Q \subseteq P$ of $p$ and $r > 0$ such that $B := B_r^F(x) \subseteq U$ and the following holds:

(a) $f_q(B)$ is open in $E$, for each $q \in Q$, and $\phi_q: B \to f_q(B)$, $\phi_q(y) := f_q(y) = f(q, y)$ is an $SC^k$-diffeomorphism.

(b) For all $q \in Q$, $y \in B$, and $s \in {]0, r - \|y - x\|]}$, we have
$$f_q(y) + A.B_{as}(0) \subseteq f_q(B_s(y)) \subseteq f_q(y) + A.B_{bs}(0).$$

(c) $W := \bigcup_{q \in Q}(\{q\} \times f_q(B))$ is open in $Z \times F$, and $\psi: W \to B$, $\psi(q, v) := \phi_q^{-1}(v)$ an $SC^k$-map. Furthermore, the map $\theta: Q \times B \to W$, $\theta(q, y) := (q, f(q, y))$ is an $SC^k$-diffeomorphism, with inverse given by $\theta^{-1}(q, v) = (q, \psi(q, v))$.

(d) $Q \times (f(p, x) + A.B_\delta(0)) \subseteq W$ for some $\delta > 0$.

In particular, for each $q \in Q$ there is a unique element $\beta(q) \in B$ such that $f(q, \beta(q)) = f(p, x)$, and the mapping $\beta: Q \to B$ so obtained is of class $SC^k$.

If $(\mathbb{K}, |.|)$ is ultrametric here and $(E, \|.\|)$ an ultrametric Banach space, then $r$ can be chosen such that (a)–(d) can be replaced with the following stronger assertions:

(a)′ $f_q(B) = f(p, x) + A.B_r(0) =: V$, for each $q \in Q$, and $\phi_q: B \to V$, $\phi_q(y) := f(q, y)$ is an $SC^k$-diffeomorphism.

(b)′ $\phi_q(B_s(y)) = \phi_q(y) + A.B_s(0)$ for all $q \in Q$, $y \in B$ and $s \in {]0, r]}$.

(c)′ $\psi: W := Q \times V \to B$, $\psi(q, v) := \phi_q^{-1}(v)$ is an $SC^k$-map, and $\theta: Q \times B \to Q \times V = W$, $\theta(q, y) := (q, f(q, y))$ is an $SC^k$-diffeomorphism, with inverse given by $\theta^{-1}(q, v) = (q, \psi(q, v))$.

**Proof of Theorems 7.3 and 7.4.** We proceed in various steps.

**7.5** *If the assertion of the Inverse Function Theorem for $SC^k$-maps is correct for some $k \in \mathbb{N} \cup \{\infty\}$, then also the Implicit Function Theorem holds for $SC^k$-maps.*

In fact, suppose that $Z$, $F$, $P$, $U$, $(p, x)$, $E$, $0 < a < 1 < b$, an $SC^k$-map $f: P \times U \to F$, and $A$ are given as described in Theorem 7.4. Define $c := \min\left\{\frac{b-1}{\|A^{-1}\|}, \frac{1-a}{\|A^{-1}\|}\right\} < \frac{1}{\|A^{-1}\|}$. Then $1 - c\|A^{-1}\| \geq a$ and $1 + c\|A^{-1}\| \leq b$. Since $f$ is strictly differentiable at $(p, x)$, there exists $r > 0$ such that $B_r^E(p, x) = B_r^Z(p) \times B_r^F(x) \subseteq P \times U$,

$$\frac{\|f(q_2, y_2) - f(q_1, y_1) - f'(p, x).(q_2 - q_1, y_2 - y_1)\|}{\|(q_2 - q_1, y_2 - y_1)\|} \leq c \text{ for all } (q_1, y_1) \neq (q_2, y_2) \in B_r^E(p, x),$$

and $\|f'_q(y) - A\| < c$ for all $(q,y) \in B^f_r(p,x)$, since $f'_q(y) = f'(q,y)(0,\bullet)$ depends continuously on $(q,y)$ (cf. Lemma 3.2 and proof of Lemma 2.2). Abbreviate $Q := B^Z_r(p)$ and $B := B^F_r(x)$. Then $f'_q(y) \in \mathrm{GL}(F)$ for all $q \in Q$ and $y \in B$. In fact, $\|A^{-1} f'_q(y) - \mathbf{1}\| \leq \|A^{-1}\| \|f'_q(y) - A\| < \|A^{-1}\| c < 1$, whence $A^{-1} f'_q(y)$ is invertible and hence so is $f'_q(y)$. Then

$$\frac{\|f_q(z) - f_q(y) - A.(z-y)\|}{\|z - y\|} = \frac{\|f(q,z) - f(q,y) - f'(p,x).(0,z-y)\|}{\|(0, z-y)\|} \leq c \qquad (27)$$

for all $q \in Q$ and $y \neq z \in B$, where $c < \frac{1}{\|A^{-1}\|}$. Thus Lemma 6.2 applies to $f|_{Q \times B}$, showing that $f_q(B)$ is open in $E$ and $\phi_q := f_q|_B^{f_q(B)}$ a homeomorphism onto its image, for each $q \in Q$; the set $W := \bigcup_{q \in Q} \{q\} \times f_q(B)$ is open in $P \times F$ and hence in $Z \times F$, and $\psi \colon W \to B$, $\psi(q,z) := \phi_q^{-1}(z)$ is continuous; the map $\theta \colon Q \times B \to W$, $\theta(q,y) := (q, f(q,y))$ is a homeomorphism, with inverse given by $\theta^{-1}(q,z) = (q, \psi(q,z))$. Furthermore, in view of (27), Lemma 6.1 applies to $f_q|_B$ for all $q \in Q$, whence (b) holds. By the $SC^k$-case of the Inverse Function Theorem, $\phi_q \colon B \to f_q(B)$ is an $SC^k$-diffeomorphism, for all $q \in Q$. Thus (a) holds. To complete the proof of (c), note that the homeomorphism $\theta \colon Q \times B \to W$ is an $SC^k$-map, whose differential $\theta'(q,y)$ at any given point $(q,y) \in Q \times B$ can be interpreted as an upper triangular $2 \times 2$-block matrix with $\mathrm{id}_Z$ and $f'_q(y)$ on the diagonal, entailing that $\theta'(q,y)$ is invertible. Hence, by the Inverse Function Theorem for $SC^k$-maps, $\theta$ restricts to an $SC^k$-diffeomorphism (onto the image) on some open neighbourhood of $(q,y)$, entailing that $\theta^{-1}$ is an $SC^k$-map on some open neighbourhood of $\theta(q,y)$. Thus $\theta$ is an $SC^k$-diffeomorphism. Since $\theta^{-1}(q,z) = (q, \psi(q,z))$ for all $(q,z) \in W$, we readily deduce that $\psi$ is an $SC^k$-map, thus completing the proof of (c).

(d) is easily established: we set $\delta := \frac{ar}{2}$. After shrinking $Q$, we may assume that $\|f(q,x) - f(p,x)\| < \delta$ for all $q \in Q$. Then, using (b) with $y := x$ and $s := r$, we see that $\{q\} \times f_q(B) \supseteq \{q\} \times (f_q(x) + A.B_{ar}(0)) \supseteq \{q\} \times (f_p(x) + A.B_\delta(0))$, for all $q \in Q$. Thus (d) holds. The assertion concerning $\beta$ is then obvious.

In the special case where $\mathbb{K}$ is an ultrametric field and $F$ an ultrametric Banach space, we establish (a)–(c) as just described, choosing however $Q$ so small that $\|f(q,x) - f(p,x)\| < r$ for all $q \in Q$. Then $f_q(B) = f_q(y) + A.B_r(0) = f_p(y) + A.B_r(0) =: V$ for all $q \in Q$, and hence $W = Q \times V$. In fact, in view of (27), we can apply Lemma 6.1 (b) to the mapping $f_q|_B$, and then use the fact that $B$ is an additive subgroup of $F$. Furthermore, again by Lemma 6.1 (b), $f_q(B_s(y)) = f_q(y) + A.B_s(0)$ for all $q \in Q$, $y \in B$, and $s \in ]0,r]$. Thus (a)′–(c)′ hold.

**7.6** *Suppose that the assertion of the Inverse Function Theorem for $SC^k$-maps is valid for all $k \in \mathbb{N}$. Then it is also valid for $k = \infty$.*

In fact, let $E$, $U$, $x \in U$, and an $SC^\infty$-map $f \colon U \to E$ be given as described in Theorem 7.3. By the $SC^1$-version of Theorem 7.3, we find $r > 0$ such that $B := B_r(x) \subseteq U$, $f(B)$ is open in $E$ and such that $f|_B^{f(B)}$ is an invertible $SC^1$-map, with inverse $g := (f|_B^{f(B)})^{-1} \colon f(B) \to B$ of class $SC^1$. Then $g$ is of class $SC^\infty$. In fact, let $k \in \mathbb{N}$. Given any $y' \in f(B)$, set $y := g(y') \in B$. Then $f'(y)$ is invertible and thus, by the $SC^k$-case of Theorem 7.3, there exists a neighbourhood $V \subseteq B$ of $y$ such that $(f|_V^{f(V)})^{-1} = g|_{f(V)}^V$ is of class $SC^k$, where $f(V)$ is a neighbourhood of $f(y) = y'$. Thus $g$ is locally $SC^k$ and thus an $SC^k$-map. As $k \in \mathbb{N}$ was arbitrary, $g$ is an $SC^\infty$-map.

**7.7** In view of **7.5** and **7.6**, in oder to establish Theorem 7.3 and Theorem 7.4, it suffices to establish the $SC^k$-case of Theorem 7.3 for all $k \in \mathbb{N}$. This we accomplish by induction. For $k = 1$, the assertion of Theorem 7.3 is covered by Proposition 7.1.

**7.8 (Induction step).** Suppose that $2 \leq k \in \mathbb{N}$ is given, and suppose that the $SC^{k-1}$-case of the Inverse Function Theorem holds. Let $E$, $U$, $x$ and an $SC^k$-map $f \colon U \to E$ be as described in Theorem 7.3. Let $r > 0$ and $B := B_r(x)$ be as described in the $SC^{k-1}$-case of the theorem. Thus $f(B)$ is open in $E$ and $f|_B^{f(B)}$ is an invertible $SC^k$-map, whose inverse $g := (f|_B^{f(B)})^{-1}$ is an $SC^{k-1}$-map. After replacing $f$ with $f|_B$, we may assume without loss of generality that $U = B$. Set $V := f(U)$. Since $g \colon V \to U$ is an $SC^{k-1}$-map, it is clear that $g^{[1]} \colon V^{[1]} \to E$ is $SC^{k-1}$ on the open subset $\{(y, z, t) \in V^{[1]} : t \neq 0\}$ of $V^{[1]}$. Thus $g^{[1]}$ will be an $SC^{k-1}$-map (and thus $g$ an $SC^k$-map) if we can show that, for every $(y_0, z_0) \in V \times E$, the mapping $g^{[1]}$ is $SC^{k-1}$ on some open neighbourhood of $(y_0, z_0, 0)$ in $V^{[1]}$. To this end, we observe first that $f \circ g = \mathrm{id}_V$ and the Chain Rule entail that $f^{[1]}(g(y), g^{[1]}(y, z, t), t) = z$ for all $(y, z, t) \in V^{[1]}$. There are open neighbourhoods $W_1 \subseteq E$ of $g(y_0)$, $W_2 \subseteq E$ of $g^{[1]}(y_0, z_0, 0)$, and $W_3 \subseteq \mathbb{K}$ of $0$ such that $W_1 \times W_2 \times W_3 \subseteq U^{[1]}$. Next, we find an open neighbourhood $P = P_1 \times P_2 \times P_3 \subseteq V^{[1]}$ of $(y_0, z_0, 0)$ such that $g(P_1) \subseteq W_1$, $g^{[1]}(P) \subseteq W_2$, and $P_3 \subseteq W_3$. By the preceding, the $SC^{k-2}$-map $\beta := g^{[1]}|_P^{W_2} \colon P \to W_2$ satisfies
$$h((y, z, t), \beta(y, z, t)) = 0 \quad \text{for all } (y, z, t) \in P, \tag{28}$$
where $h \colon P \times W_2 \to E$ is the $SC^{k-1}$-map defined via $h((y, z, t), w) := f^{[1]}(g(y), w, t) - z$. Since $h((y_0, z_0, 0), w) = f'(g(y_0)).w - z_0$ is affine-linear in $w$, the differential of $h$ with respect to the $w$-variable satisfies $A := d_2 h((y_0, z_0, 0), \beta(y_0, z_0, 0); \bullet) = f'(g(y_0)) \in \mathrm{GL}(E)$. Since $\beta$ is a continuous solution to the implicit equation (28), we deduce from the $SC^{k-1}$-case of the Implicit Function Theorem 7.4 (which holds in view of the induction hypothesis and **7.5**) that $\beta$ is $SC^{k-1}$ on some open neighbourhood of $(y_0, z_0, 0)$, as we set out to show. This completes the proof of Theorems 7.3 and 7.4. $\square$

**Remark 7.9** If $k > 1$, in the preceding induction step we encounter a map $\beta$ defined on $P \subseteq V^{[1]} \subseteq E^{[1]}$. In order that $E^{[1]} = E \times E \times \mathbb{K}$ be a Banach space (so that the implicit function theorem can be applied), it is necessary that $\mathbb{K}$ be complete.

# 8 General Ultrametric Implicit Function Theorem

We are now in the position to prove a generalized implicit function theorem for mappings from open subsets of metrizable topological vector spaces over complete ultrametric fields to Banach spaces over such fields.

**Theorem 8.1 (Ultrametric Implicit Function Theorem)** *Let $(\mathbb{K}, |.|)$ be a complete ultrametric field, $k \in \mathbb{N} \cup \{\infty\}$, $Z$ be a metrizable topological $\mathbb{K}$-vector space, and $E$ be a Banach space over $\mathbb{K}$. Let $P \subseteq Z$ and $U \subseteq E$ be open subsets, and $f\colon P \times U \to E$ be a map. We assume that at least one of the following conditions is satisfied:*

(i) *$\mathbb{K}$ is locally compact, $E$ is finite-dimensional, and $f$ is of class $C_\mathbb{K}^k$. Or:*

(ii) *$f$ is of class $C_\mathbb{K}^{k+1}$.*

*We abbreviate $f_q := f(q, \bullet)\colon U \to E$ for $q \in P$. Suppose that $(p, x) \in P \times U$ is given such that $A := f_p'(x) := d_2 f(p, x, \bullet) := df((p, x), (0, \bullet)) \in \mathrm{GL}(E)$. Let $a, b \in \mathbb{R}$ be given such that $0 < a < 1 < b$. Then there exists an open neighbourhood $Q \subseteq P$ of $p$ and $r > 0$ such that $B := B_r(x) \subseteq U$ and the following holds:*

(a) *$f_q(B)$ is open in $E$, for each $q \in Q$, and $\phi_q\colon B \to f_q(B)$, $\phi_q(y) := f_q(y) = f(q, y)$ is an $SC^k$-diffeomorphism.*

(b) *For all $q \in Q$, $y \in B$, and $s \in \,]0, r - \|y - x\|]$, we have*
$$f_q(y) + A.B_{as}(0) \subseteq f_q(B_s(y)) \subseteq f_q(y) + A.B_{bs}(0).$$

(c) *$W := \bigcup_{q \in Q}(\{q\} \times f_q(B))$ is open in $Z \times E$, and the map $\psi\colon W \to B$, $\psi(q, v) := \phi_q^{-1}(v)$ is of class $C^k$. Furthermore, the map $\theta\colon Q \times B \to W$, $\theta(q, y) := (q, f(q, y))$ is a $C_\mathbb{K}^k$-diffeomorphism, with inverse given by $\theta^{-1}(q, v) = (q, \psi(q, v))$.*

(d) *$Q \times (f_p(x) + A.B_\delta(0)) \subseteq W$ for some $\delta > 0$.*

*In particular, for each $q \in Q$ there is a unique element $\beta(q) \in B$ such that $f(q, \beta(q)) = f(p, x)$, and the mapping $\beta\colon Q \to B$ so obtained is of class $C_\mathbb{K}^k$.*

*If $(E, \|.\|)$ is an ultrametric Banach space here, then $Q$ and $r$ can be chosen such that (a)–(d) can be replaced with the following stronger assertions:*

(a)′ *$f_q(B) = f(p, x) + A.B_r(0) =: V$, for each $q \in Q$, and $\phi_q\colon B \to V$, $\phi_q(y) := f(q, y)$ is an $SC^k$-diffeomorphism.*

(b)′ *$f_q(B_s(y)) = f_q(y) + A.B_s(0)$ for all $q \in Q$, $y \in B$ and $s \in \,]0, r]$.*

(c)′ *The mapping $\psi\colon Q \times V \to B$, $\psi(q, v) := \phi_q^{-1}(v)$ is of class $C_\mathbb{K}^k$. Furthermore, $\theta\colon Q \times B \to Q \times V$, $\theta(q, y) := (q, f(q, y))$ is a $C_\mathbb{K}^k$-diffeomorphism, with inverse given by $\theta^{-1}(q, v) = (q, \psi(q, v))$.*

**Proof.** Define $c := \min\left\{\frac{b-1}{\|A^{-1}\|}, \frac{1-a}{\|A^{-1}\|}\right\} < \frac{1}{\|A^{-1}\|}$, where $A := f_p'(x)$. Then $1 - c\|A^{-1}\| \geq a$ and $1 + c\|A^{-1}\| \leq b$. In the situation of (i), let $e_1, \ldots, e_n$ be a basis of $E$. The mappings $P \times U \to E$, $(q, y) \mapsto d_2 f(q, y, e_i)$ being continuous for $i = 1, \ldots, n$, also the map

$$P \times U \to L(E), \quad (q, y) \mapsto f_q'(y) = d_2 f(q, y, \bullet) \tag{29}$$

is continuous. In the situation of (ii), the map in (29) is continuous as well, by Lemma 3.5. In either case, since $\mathrm{GL}(E)$ is open in $L(E)$ and $f'_p(x) \in \mathrm{GL}(E)$, after replacing $P$ and $U$ with smaller open neighbourhoods of $p$ and $x$, respectively, we may assume that $f'_q(y) \in \mathrm{GL}(E)$ for all $(q,y) \in P \times U$, and $\|f'_q(y) - f'_p(x)\| \leq \frac{c}{2}$. Using Lemma 4.5 (resp., Lemma 3.5), we find an open neighbourhood $Q \subseteq P$ of $p$ and $r > 0$ such that $B := B_r(x) \subseteq U$ and

$$\frac{\|f_q(z) - f_q(y) - f'_q(x).(z-y)\|}{\|z-y\|} \leq \frac{c}{2} \tag{30}$$

for all $y \neq z \in B$. As a consequence,

$$\frac{\|f_q(z) - f_q(y) - f'_p(x).(z-y)\|}{\|z-y\|} \leq \frac{\|f_q(z) - f_q(y) - f'_q(x).(z-y)\|}{\|z-y\|} + \|f'_p(x) - f'_q(x)\|$$
$$\leq c$$

for all $q \in Q$ and $z \neq y \in B$, entailing that

$$\sup\left\{\frac{\|f_q(z) - f_q(y) - f'_p(x).(z-y)\|}{\|z-y\|} : q \in Q, z \neq y \in B\right\} \leq c < \frac{1}{\|A^{-1}\|}. \tag{31}$$

Thus Lemma 6.2 applies to $f|_{Q \times B}$, showing that $f_q(B)$ is open in $E$ and $\phi_q := f_q|_B^{f_q(B)}$ a homeomorphism onto its image, for each $q \in Q$; the set $W := \bigcup_{q \in Q}\{q\} \times f_q(B)$ is open in $P \times E$ and thus in $Z \times E$, and $\psi : W \to B$, $\psi(q,z) := \phi_q^{-1}(z)$ is continuous; the map $\theta : Q \times B \to W$, $\theta(q,y) := (q, f(q,y))$ is a homeomorphism, with inverse given by $\theta^{-1}(q,z) = (q, \psi(q,z))$. Furthermore, in view of (31), Lemma 6.1 applies to $f_q|_B$, for all $q \in Q$, showing that (b) holds. By the Inverse Function Theorem for $SC^k$-maps, the map $\phi_q : B \to f_q(B)$ is an $SC^k$-diffeomorphism, for all $q \in Q$. Thus (a) holds.

(d) is easily established: we set $\delta := \frac{ar}{2}$. After shrinking $Q$, we may assume that $\|f(q,x) - f(p,x)\| < \delta$ for all $q \in Q$. Then, using (b) with $y := x$ and $s := r$, we get $\{q\} \times f_q(B) \supseteq \{q\} \times (f_q(x) + A.B_{ar}(0)) \supseteq \{q\} \times (f_p(x) + A.B_\delta(0))$, for all $q \in Q$. Thus (d) holds.

To see that the continuous map $\psi$ is of class $C^k_\mathbb{K}$ (which will entail the validity of (c)), in view of Proposition 1.15, it suffices to show that, for all smooth maps $c : \mathbb{K}^{k+1} \to W$, the composition $\psi \circ c : \mathbb{K}^{k+1} \to E$ is of class $C^k_\mathbb{K}$. Since $W \subseteq Q \times E$, we have $c = (c_1, c_2)$ with smooth mappings $c_1 : \mathbb{K}^{k+1} \to Q \subseteq Z$, $c_2 : \mathbb{K}^{k+1} \to E$. Define $h : \mathbb{K}^{k+1} \times B \to V$, $h(t,y) := f(c_1(t), y)$. Then $h$ is a $C^k_\mathbb{K}$-map in the situation of (i) and hence $SC^k$ (Remark 5.4). In the situation of (ii), $h$ is of class $C^{k+1}_\mathbb{K}$ and hence $SC^k$ (see Remark 5.3). Given $t \in \mathbb{K}^{k+1}$, abbreviate $h_t := h(t, \bullet) : B \to E$; by the above, $h_t$ has open image and is a homeomorphism onto its image. Since $h'_t(y) = d_2f(c_1(t), y, \bullet) \in \mathrm{GL}(E)$ for all $(t,y) \in \mathbb{K}^{k+1} \times B$, we deduce from the Implicit Function Theorem for $SC^k$-maps (Theorem 7.4) that $\kappa : W \to B$, $\kappa(t,z) := h_t^{-1}(z)$ is an $SC^k$-map and hence of class $C^k_\mathbb{K}$. Now $\psi(c_1(t), c_2(t)) = \kappa(t, c_2(t))$ for all $t \in \mathbb{K}^{k+1}$ shows that $\psi \circ c$ is of class $C^k_\mathbb{K}$, as required. Thus $\psi$ is $C^k_\mathbb{K}$.

In the special case where $E$ is an ultrametric Banach space, we establish (a)–(c) as just described, choosing however $Q$ so small that $\|f(q,x) - f(p,x)\| < r$ for all $q \in Q$ (we might actually replace $\frac{c}{2}$ with $c$ in (30) now). Then $f_q(B) = f_q(y) + A.B_r(0) = f_p(y) + A.B_r(0) =:$
$V$ for all $q \in Q$ (by Lemma 6.1 (b), applied as at the end of **7.5**), and hence $W = Q \times V$.
Furthermore, again by Lemma 6.1 (b), $f_q(B_s(y)) = f_q(y) + A.B_s(0)$ for all $q \in Q$, $y \in B$,
and $s \in \,]0,r]$. Thus (a)′–(c)′ hold. □

**Remark 8.2** Three cases described in the table given in the introduction still remain to be discussed.

(a) Suppose we retain the hypotheses of Theorem 8.1, with $k = 1$, except that we let $Z$ be an arbitrary topological $\mathbb{K}$-vector space now (which need not be metrizable). Suppose we are in the situation of (i). Then the proof of Theorem 8.1 shows that the following weakened conclusions of the theorem remain valid: (a), (b) and (d) will hold unchanged; $\psi$ in (c) and $\beta$ will be continuous; $\theta$ in (c) will be a homeomorphism (likewise for (a)′–(c)′).

(b) Suppose we retain the hypotheses of Theorem 8.1, with $k = 1$, except that we let $\mathbb{K}$ be an arbitrary valued field and $Z$ be an arbitrary topological $\mathbb{K}$-vector space. Suppose we are in the situation of (ii). Then the proof of Theorem 8.1 shows that the following weakened conclusions remain valid: (a), (b) and (d) will hold unchanged; $\psi$ in (c) and $\beta$ will be continuous; $\theta$ in (c) will be a homeomorphism.

(c) Suppose we retain the hypotheses of Theorem 8.1, with $k = 1$, except that we let $\mathbb{K}$ be a subfield of $\mathbb{R}$ now, equipped with the absolute value obtained by restricting the usual absolute value on $\mathbb{R}$. Suppose we are in the situation of (ii). Then (a)–(d) and their proof remain valid verbatim, and $\beta$ is $C^1$. In fact, Proposition 1.15 remains valid when $\mathbb{R}$ is replaced with arbitrary subfields of $\mathbb{R}$ (the proof given in [2] applies without changes).

# Appendix: $FC^k$-maps vs. $SC^k$-maps in the real case

**Theorem A.3** *Suppose that $E$ is a normed vector space over $\mathbb{R}$, $F$ a real locally convex topological vector space, $U \subseteq E$ an open subset, $f \colon U \to F$ a mapping, and $k \in \mathbb{N}_0 \cup \{\infty\}$. If $f$ is $FC^k$, then $f$ is $SC^k$.*

**Proof.** We may assume that $k \in \mathbb{N}_0$. The proof is by induction. The case $k = 0$ is trivial, and the case $k = 1$ is a standard fact (see [5, 2.3.3], cf. also [7, Thm. 3.8.1]).

*Induction step.* Suppose that $k \geq 2$, and suppose that every $FC^{k-1}$-map is $SC^{k-1}$. Let $f \colon E \supseteq U \to F$ be an $FC^k$-map. Then $f$ is $SC^{k-1}$ and hence $SC^1$ in particular. Then, $f$ being $SC^{k-1}$, so is $f^{[1]}$ on $\{(x,y,t) \in U^{[1]} : t \neq 0\}$. It therefore only remains to show that, for every $x_0 \in U$ and $y_0 \in E$, the map $f^{[1]}$ is $SC^{k-1}$ on some open neighbourhood of $(x_0, y_0, 0)$. There is $r > 0$ such that $B_{2r}(x_0) \subseteq U$. Choose $\delta \in ]0, r[$ such that $\delta < \frac{r}{2(\|y_0\|+1)}$. Then $V := B_r(x_0) \times B_1(y_0) \times ]-\delta, \delta[ \subseteq U^{[1]}$, and we have

$$f^{[1]}(x, y, t) = \int_0^1 df(x + sty, y)\, ds = \int_0^1 h((x, y, t), s)\, ds \quad \text{for all } (x, y, t) \in V, \qquad (32)$$

where $h \colon V \times ]-2\delta, 2\delta[ \to F$, $h((x, y, t), s) := df(x + sty, y)$ is an $FC^{k-1}$-map. In view of (32), we inductively deduce from [8, 8.11.2] that $f^{[1]}|_V$ is an $FC^{k-1}$-map, if $F$ is a *Banach space*. If, more generally, $F$ is a *complete* locally convex space, then $F = \varprojlim F_i$ is a projective limit in the category of locally convex spaces of some projective system of Banach spaces, and thus $L(H, F) = \varprojlim L(H, F_i)$ (which again is a complete locally convex space), for every normed space $H$. A simple inductive argument now shows that a map $g$ from an open subset of a normed space to a projective limit $F = \varprojlim F_i$ is an $FC^k$-map if and only if $\pi_i \circ g$ is $FC^k$ for each $i$, where $\pi_i \colon F \to F_i$ are the limit maps. In the situation we are interested in, $\pi_i \circ f^{[1]}|_V = (\pi_i \circ f)^{[1]}|_V$ maps into a Banach space and hence is an $FC^{k-1}$-map, by what has already been shown, whence $f^{[1]}|_V$ is an $FC^{k-1}$-map to the projective limit $F$. In the *general case*, when $F$ is not necessarily complete, the preceding shows that $f^{[1]}|_V$ is $FC^{k-1}$ as a mapping into the completion $\widetilde{F}$ of $F$. Since $(f^{[1]})'(x) = d(f^{[1]})(x, \bullet)$ for $x \in V$ actually is a mapping into $F$ (not only into $\widetilde{F}$), and likewise for the higher order differentials, we deduce that $f^{[1]}|_V$ is $FC^{k-1}$ as a mapping into $F$ also in the fully general case.

Now $f^{[1]}|_V$ being an $FC^{k-1}$-map, it is an $SC^{k-1}$-map, by the induction hypothesis. Thus $f$ is $SC^1$ with $f^{[1]}$ an $SC^{k-1}$-map, and hence $f$ is $SC^k$. $\square$

The author does not know whether, conversely, every $SC^k$-map is $FC^k$. For $k = 1$, this is well-known [5, 2.3.3], but the generalization to higher $k$ does not seem to be clear.

[24] Leslie J (1982) On the group of real analytic diffeomorphisms of a compact real analytic manifold. Trans Amer Math Soc **274**: 651–669

[25] Losik M V (1992) Fréchet manifolds as diffeologic spaces. Russ Math **36** no 5: 31–37

[26] Ludkovsky S V (1996) Measures on groups of diffeomorphisms of non-archimedian Banach manifolds. Russian Math Surv **51** no 2: 338–340

[27] —— (2000) Quasi-invariant measures on non-Archimedian groups and semigroups of loops and paths, their representations I. Ann Math Blaise Pascal **7** no 2: 19–53

[28] Ma T W (2001) Inverse mapping theorem on coordinate spaces. Bull London Math Soc **33**: 473–482

[29] Monna A F (1979) Analyse Non-Archimédienne. Springer

[30] Rooij A C M (1978) Non-Archimedian Functional Analysis. Marcel Dekker

[31] Schikhof W H (1970) Differentiation in non-archimedian valued fields. Nederl Akad Wet Proc Ser A **73**: 35–44

[32] —— (1984) Ultrametric Calculus. Cambridge University Press

[33] Serre J-P (1992) Lie Algebras and Lie Groups. Springer

[34] Slezák B (1988) On the inverse function theorem and implicit function theorem in Banach spaces. In: Musielak J (Ed) Function Spaces (Poznań, 1986), pp 186–190. Leipzig: Teubner

[35] De Smedt S (1998) Local invertibility of non-archimedian vector-valued functions. Ann Math Blaise Pascal **5** no 1: 13–23

[36] Souriau J-M (1984) Groupes différentiels de physique mathématique. In: Dazord P, Desolneux-Moulis N (Eds) Feuilletages et quantification géometrique. Journ lyonnaises Soc math France 1983, Sémin sud-rhodanien de géom II, pp 73–119. Paris: Hermann

[37] Weil A (1973) Basic Number Theory. Springer

[38] Więsław W (1988) Topological Fields. New York and Basel: Marcel Dekker
**Helge Glöckner**, TU Darmstadt, FB Mathematik AG 5, Schlossgartenstr. 7, 64289 Darmstadt, Germany. E-Mail: gloeckner@mathematik.tu-darmstadt.de
36